# Graphical Abstract

## MPMICE: A hybrid MPM-CFD model for simulating coupled problems in porous media. Application to earthquake-induced submarine landslides

Quoc Anh Tran, Gustav Grimstad, Seyed Ali Ghoreishian Amiri

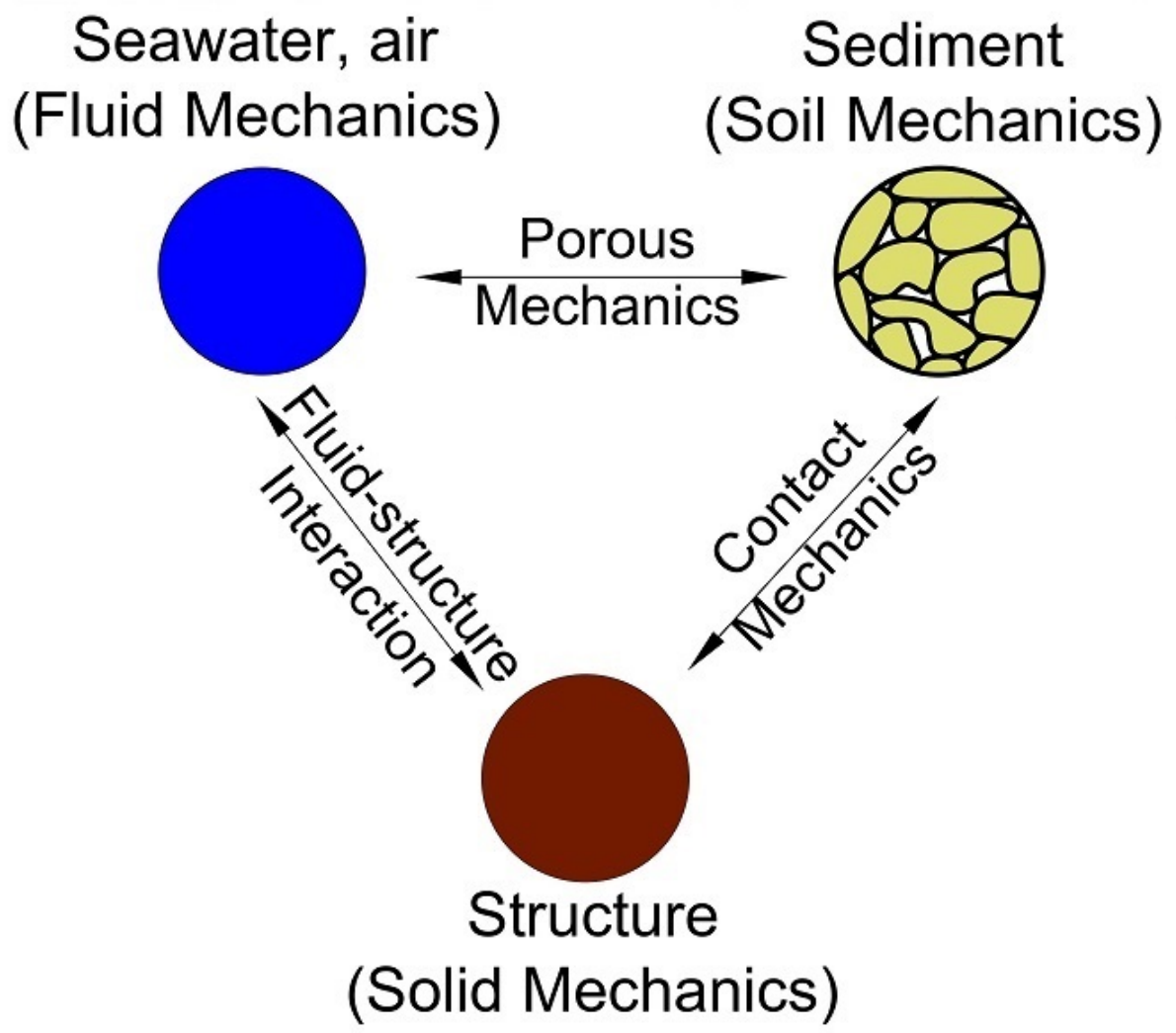

**soil-fluid-structure interaction**

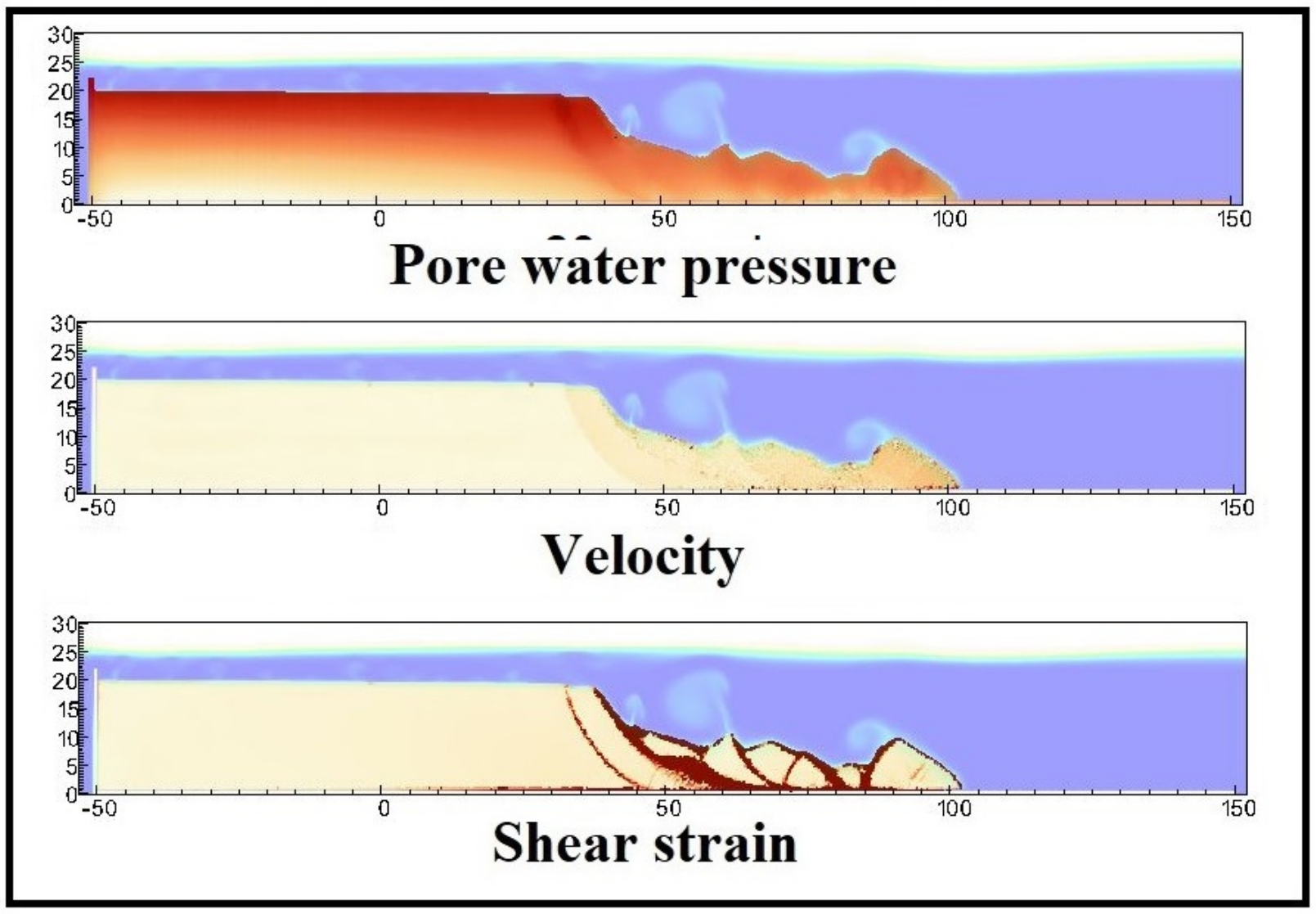

**Application to earthquake-induced submarine landslide**

# Highlights

**MPMICE: A hybrid MPM-CFD model for simulating coupled problems in porous media. Application to earthquake-induced submarine landslides**

Quoc Anh Tran, Gustav Grimstad, Seyed Ali Ghoreishian Amiri

- MPMICE is introduced for multiphase flow in porous media.
- Material Point method allows to model large deformation of non-isothermal porous media.
- ICE (compressible multi-material CFD formulation) allows to stablize pore water pressure and turbulent flow.
- MPMICE is validated and apply to simulate the earthquake-induced submarine landslide.

# MPMICE: A hybrid MPM-CFD model for simulating coupled problems in porous media. Application to earthquake-induced submarine landslides

Quoc Anh Tran[a], Gustav Grimstad[a], Seyed Ali Ghoreishian Amiri[a]

[a]*Norwegian University of Science and Technology, , Trondheim, 7034, Norway*

---

**Abstract**

In this paper, we describe a soil-fluid-structure interaction model that combines soil mechanics (saturated sediments), fluid mechanics (seawater or air), and solid mechanics (structures). The formulation combines the Material Point Method, which models large deformation of the porous media and the structure, with the Implicit Continuous-fluid Eulerian, which models complex fluid flows. We validate the model and simulate the whole process of earthquake-induced submarine landslides. We show that this model captures complex interactions between saturated sediment, seawater, and structure, so we can use the model to estimate the impact of potential submarine landslides on offshore structures.

---

Contents

## Nomenclature

**General variables**

| Variable | Dimensions | Description |
|---|---|---|
| $V$ | $[L^3]$ | Representative volume |
| $n$ | | Porosity |
| $\boldsymbol{\sigma}$ | $[ML^{-1}t^{-2}]$ | Total stress tensor |
| $\Delta t$ | $[t]$ | Time increment |
| $\boldsymbol{b}$ | $[ML^1t^{-2}]$ | Body force |
| $c_v$ | $[L^2t^{-2}T^{-1}]$ | Constant volume specific heat |
| $f_d$ | $[MLt^{-2}]$ | Drag forces in momentum exchange term |
| $f^{int}$ | $[MLt^{-2}]$ | Internal forces |
| $f^{ext}$ | $[MLt^{-2}]$ | External forces |
| $q_{fs}$ | $[MLt^{-2}]$ | Heat exchange term |
| $S$ | | Weighting function |
| $\nabla S$ | | Gradient of weighting function |

**Solid phase**

| Variable | Dimensions | Description |
|---|---|---|
| $m_s$ | $[M]$ | Solid mass |
| $\rho_s$ | $[ML^{-3}]$ | Solid density |
| $\phi_s$ | | Solid volume fraction |
| $\overline{\rho}_s$ | $[ML^{-3}]$ | Bulk Solid density |
| $\boldsymbol{x}_s$ | $[L]$ | Solid Position vector |
| $\boldsymbol{U}_s$ | $[Lt^{-1}]$ | Solid Velocity vector |
| $\boldsymbol{a}_s$ | $[Lt^{-2}]$ | Solid Acceleration vector |
| $\boldsymbol{\sigma}'$ | $[ML^{-1}t^{-2}]$ | Effective Stress tensor |
| $\boldsymbol{\epsilon}_s$ | | Strain tensor |
| $e_s$ | $[L^2t^{-2}]$ | Solid Internal energy per unit mass |
| $T_s$ | $[T]$ | Solid Temperature |
| $\boldsymbol{F}_s$ | | Solid Deformation gradient |
| $V_s$ | $[L^3]$ | Solid Volume |

**Fluid phase**

| Variable | Dimensions | Description |
|---|---|---|
| $m_f$ | $[M]$ | Fluid mass |
| $\rho_f$ | $[ML^{-3}]$ | Fluid density |
| $\phi_f$ | | Fluid volume fraction |
| $\overline{\rho}_f$ | $[ML^{-3}]$ | Bulk Fluid density |
| $\boldsymbol{U}_f$ | $[Lt^{-1}]$ | Fluid Velocity vector |
| $\boldsymbol{\sigma}_f$ | $[ML^{-1}t^{-2}]$ | Fluid stress tensor |
| $p_f$ | $[ML^{-1}t^{-2}]$ | Fluid isotropic pressure |
| $\boldsymbol{\tau}_f$ | $[ML^{-1}t^{-2}]$ | Fluid shear stress tensor |
| $e_f$ | $[L^2t^{-2}]$ | Fluid Internal energy per unit mass |
| $T_f$ | $[T]$ | Fluid Temperature |
| $\upsilon_f$ | $[L^3/M]$ | Fluid Specific volume $\frac{1}{\rho_f}$ |
| $\alpha_f$ | $[1/T]$ | Thermal expansion |
| $\mu$ | $[ML^{-1}t^{-1}]$ | Fluid vicousity |
| $V_f$ | $[L^3]$ | Fluid Volume |

**Superscript**

| Variable | Dimensions | Description |
|---|---|---|
| $n$ | | Current time step |
| $L$ | | Lagrangian values |
| $n+1$ | | Next time step |

**Subscript**

| | | |
|---|---|---|
| $c$ | | Cell-centered quantity |
| $p$ | | Particle quantity |
| $i$ | | Node quantity |
| $FC$ | | Face-centered quantity |
| $L, R$ | | Left and Right cell faces |

## Introduction

Many geological natural processes and their interactions with man-made structures are influenced by soil-fluid-structure interactions. The prediction of these processes requires a tool that can capture complex interactions between soil, fluid, and structure, such as the process of submarine landslides. Indeed, The offshore infrastructure as well as coastal communities may be vulnerable to submarine landslides. Submarine landslides contain three stages: triggering, failure, and post-failure. Erosion or earthquakes can trigger slope failures in the first stage. Following the failure, sediments move quickly after the post-failure stage. In other words, solid-like sediments will behave like a fluid after failure. This phase transition makes the simulation of submarine landslides a challeging task.

Due to this phase transition, submarine landslide can be modeled by either the Computational Fluid Dynamics (CFD) or the particle-based methods. For simulating submarine slides, CFD methods solve governing equations in a full-Eulerian framework [1, 2, 3, 4] with interface capturing techniques. While CFD can handle complex flows (such as turbulent flows), it cannot account for the triggering mechanism of submarine landslides because it is not straightforwad to consider 'soil constitutive laws' of sediment materials in the Eulerian framework. In contrast, particle-based methods can overcome this problem by using the Lagrangian framework. These methods have been extensively used to simulate landslides, like Material Point Method (MPM) [5], Smooth Particle Hydro Dynamics [6], Particle Finite Element Method [7], or Coupled Eulerian Lagrangian Method [8]. For simplicity, these simulations adopt total stress analysis which neglects the pore pressure development which is key factor triggering slope failures.

Recent developments in particle-based methods model the coupling of fluid flows in porous media by sets of Lagrangian particles. For the MPM family, it is the double-point MPM ([9, 10, 11]) where fluid particles and solid particles are overlapped in a single computational grid. Even if fluid flows are considered, particle-based methods have numerical instability in modeling the fluid flow, which requires additional numerical treatments such as the B-bar method [9], null-space filter [12], or least square approximation [13, 14]. Indeed, CFD is a more optimal option for complex fluid flows especially dealing with large distortions of continuum fluid media. Therefore, it could be ideal to combine the CFD with particle-based methods. More than 50 particle-based methods have been developed to solve large deformations of solids over the last two decades [15], but the MPM appears to be the best candidate for coupling with the CFD. Because MPM incorporates a stationary mesh during computation, just like CFD. As such, both MPM and CFD can be coupled naturally in a unified computational mesh.

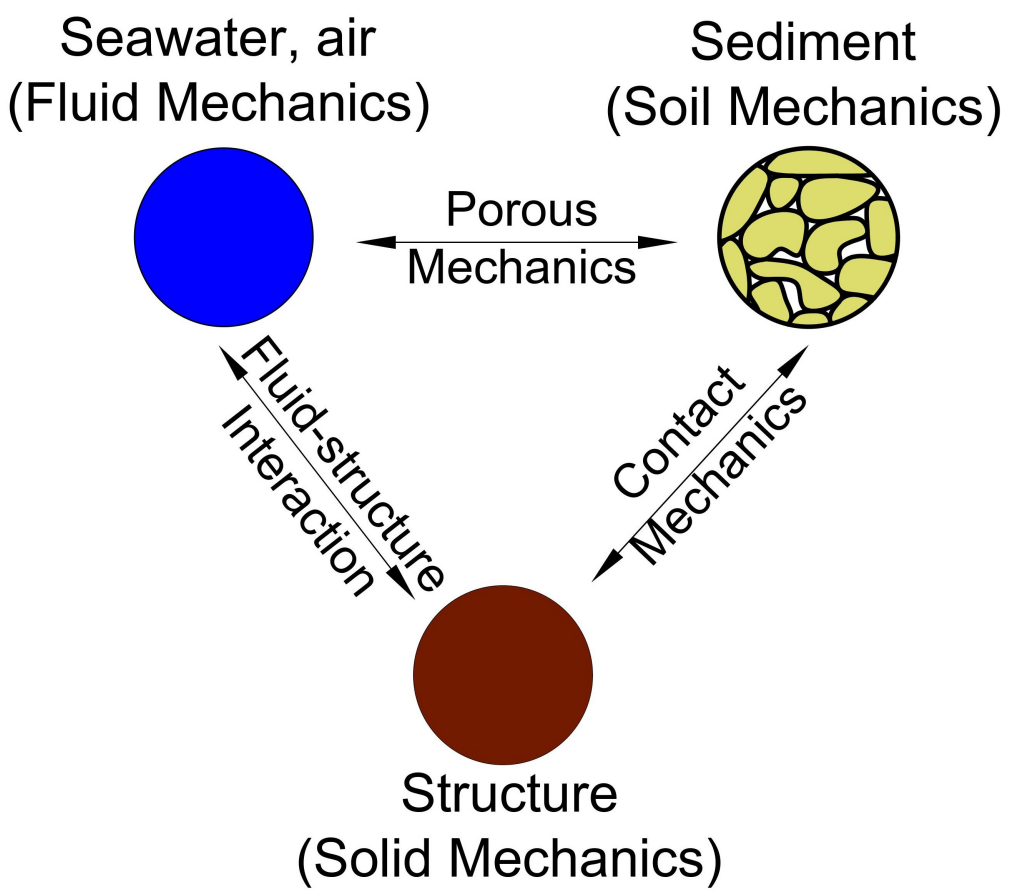

*Figure 1: Interaction between soil-fluid-structure*

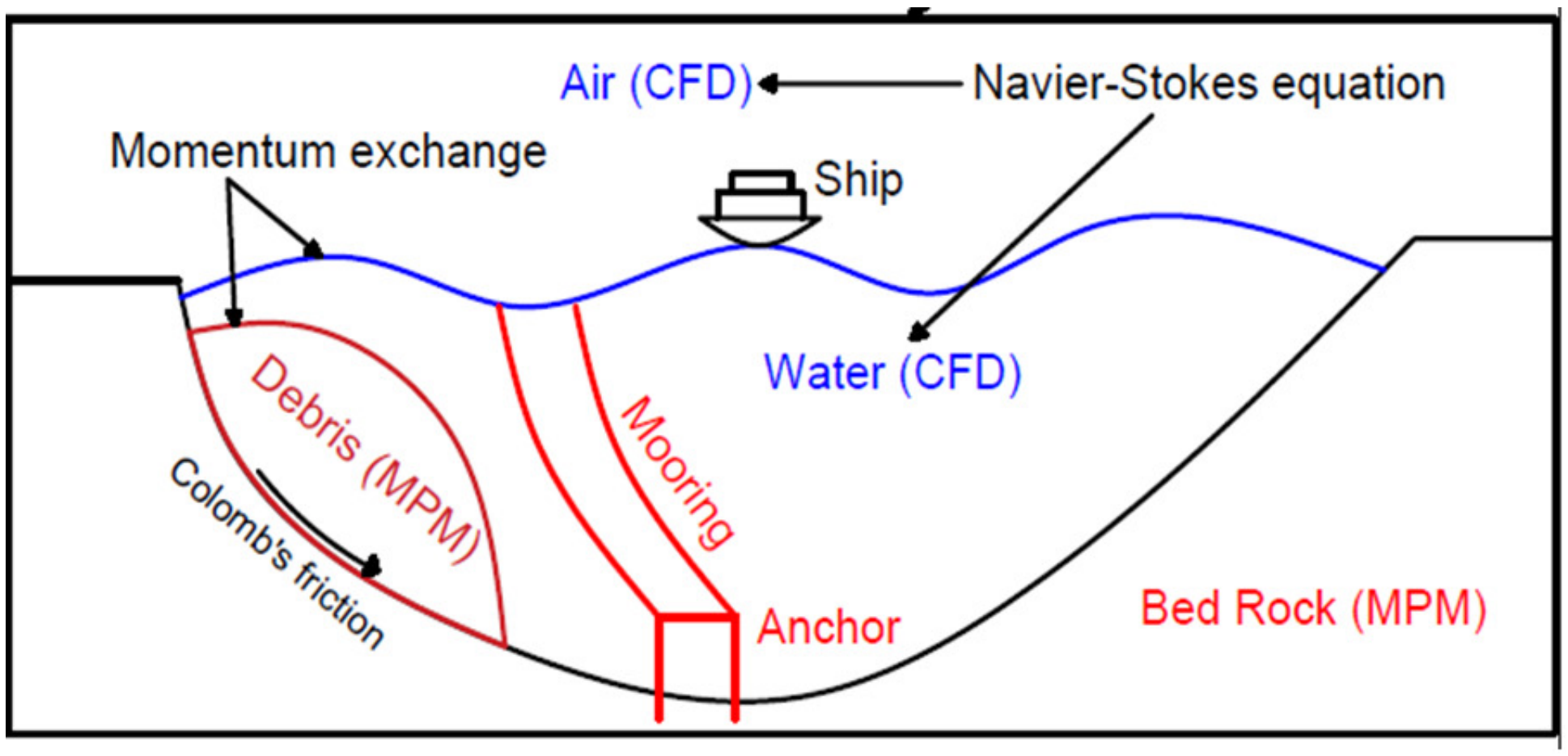

*Figure 2: Coupling of soil-water-structure interaction using MPMICE*

A numerical method for simulating soil-fluid-structure interaction (Figure 1) involving large deformations, is presented in this work in order to simulate the interaction between sediment (soil), seawater (fluid) and offshore structures (structure) namely MPMICE (Figure 2). In the MPMICE, the Material Point Method (MPM) is coupled with the Implicit Continuous Eulerian (ICE). The MPM method is a particle method that allows the porous soil to undergo arbitrary distortions. The ICE method, on the other hand, is a conservative finite volume technique with all state variables located at the cell center (temperature, velocity, mass, pressure). An initial technical report [16] at Los Alamos National Laboratory provided the theoretical and algorithmic foundation for the MPMICE, followed by the MPMICE development and implementation in the high-performance Uintah computational framework for simulating fluid-structure interactions [17]. This paper primarily contributes futher to the development of the MPMICE for analyzing the **soil**-fluid-structure interaction, since sediment should be considered as a porous media (soil) and not as a solid to capture the evolution of the pore water pressure. Baumgarten et al. [18] made the first attempt at coupling the Finite Volume Method with the MPM for the simulation of soil-fluid interaction. In contrast to the mentioned work, we use implicit time integration for the multi phase flows instead of explicit time integration for the single-phase flow.

## Theory and formulation

This section lay out the theoretical framework for the MPMICE model. We use the common notation of the continuum mechaniccs with vector and tensor denoted simply by using bold font and scalar denoted by using normal font. The notation are shown in Nomenclature.

### *Assumptions*

The following assumptions are made for the MPMICE model.

1. Solid phases (MPM) are described in a Lagrangian formulation while fluid phases (ICE) are described in an Eulerian formulation in the framework of continuum mechanics and mixture theory.
2. Solid grains are incompressible while the fluid phases are compressible. Solid's thermal expasion is negligible.
3. There is no mass exchange between solid and fluid phases.
4. Terzaghi's effective stress is valid.

*Governing equations*

A representative element volume $\Omega$ is decomposed by two domains: solid domains $\Omega_s$ and fluid domains $\Omega_f$. Then, all domains are homogenized into two overlapping continua. Considering the volume fraction of solid $\phi_s = \Omega_s/\Omega$ and fluid $\phi_f = \Omega_f/\Omega$ with the true (or Eulerian) porosity $n = \sum \phi_f$ of the representative element volume, the average density of solid and fluid phases are defined as:

$$\overline{\rho}_s = \phi_s \rho_s, \qquad \overline{\rho}_f = \phi_f \rho_f \tag{1}$$

The mass of solid and fluid phases are:

$$m_s = \int_{\Omega_s} \rho_s dV = \overline{\rho}_s V, \qquad m_f = \int_{\Omega_f} \rho_f dV = \overline{\rho}_f V \tag{2}$$

Reviewing the Terzaghi's effective stress concept for the saturated porous media, the total stress $\boldsymbol{\sigma}$ is calculated by:

$$\boldsymbol{\sigma} = \boldsymbol{\sigma'} - p_f \boldsymbol{I} \tag{3}$$

The balance equations are derived based on the mixture theory. The representative thermodynamic state of the fluid phases are given by the vector $[m_f, \boldsymbol{U}_f, e_f, T_f, v_f]$ which are mass, velocity, internal energy, temperature, specific volume. The representative state of the solid phases are given by the vector $[m_s, \boldsymbol{U}_s, e_s, T_s, \boldsymbol{\sigma'}, p_f]$ which are mass, velocity, internal energy, temperature, effective stress and pore water pressure. The derivation is presented in detail in the Appendix.

---

Mass Conservation

The mass balance equations for both fluid (e.g., water, air) and solid phases are:

$$\frac{1}{V}\frac{\partial m_f}{\partial t} + \nabla \cdot (m_f \boldsymbol{U}_f) = 0, \frac{1}{V}\frac{D_s m_s}{Dt} = 0 \tag{4}$$

Solving the mass balance equation leads to:

$$\frac{D_s n}{Dt} = \phi_s \nabla \cdot \boldsymbol{U}_s \tag{5}$$

---

Momentum Conservation

The momentum balance equation for the fluid phases (e.g., water, air) are:

$$\frac{1}{V}\left[\frac{\partial(m_f\boldsymbol{U}_f)}{\partial t} + \nabla\cdot(m_f\boldsymbol{U}_f\boldsymbol{U}_f)\right] = -\phi_f\nabla p_f + \nabla\cdot\boldsymbol{\tau}_f + \overline{\rho}_f\boldsymbol{b} + \sum\boldsymbol{f}_d \tag{6}$$

The momentum balance equation for the solid phases are:

$$\frac{1}{V}\frac{D_s(m_s\boldsymbol{U}_s)}{Dt} = \nabla\cdot(\boldsymbol{\sigma'}) - \phi_s\nabla p_f + \overline{\rho}_s\boldsymbol{b} + \sum\boldsymbol{f}_{fric} - \sum\boldsymbol{f}_d \tag{7}$$

---

Energy Conservation

The internal energy balance equation for the fluid phases (e.g., water, air) are:

$$\frac{1}{V}\left[\frac{\partial(m_f e_f)}{\partial t} + \nabla\cdot(m_f e_f\boldsymbol{U}_f)\right] = -\overline{\rho}_f p_f\frac{D_f\upsilon_f}{Dt} + \boldsymbol{\tau}_f : \nabla\boldsymbol{U}_f + \nabla\cdot\boldsymbol{q}_f + \sum q_{sf} \tag{8}$$

The internal energy balance equation for the solid phase is:

$$\frac{m_s}{V}c_v\frac{D_s(T_s)}{Dt} = \boldsymbol{\sigma'} : \frac{D_s(\boldsymbol{\epsilon}_s^p)}{Dt} + \nabla\cdot\boldsymbol{q}_s - \sum q_{sf} \tag{9}$$

where $c_v$ is the specific heat at constant volume of the solid materials.

---

Closing the systems of equations, the following additional models are needed:
(1) A constitutive equation to describe the stress - strain behaviour of solid phase (computing effective stress $\boldsymbol{\sigma'}$).
(2) Optional turbulent model to compute the viscous shear stress $\boldsymbol{\tau}_f$.
(3) Frictional forces $\boldsymbol{f}_{fric}$ for the contact for soil-structure interaction between solid/porous materials with the friction coefficient $\mu_{fric}$.
(4) Exchange momentum models (computing drag force $\boldsymbol{f}_d$) for interaction between materials.
(5) Energy exchange models (computing temerature exhange term $q_{sf}$) for interaction between materials.
(6) An equation of state to establish relations between thermodynamics variables of each fluid materials $[P_f, \overline{\rho}_f, \upsilon_f, T_f, e_f]$.

Four thermodynamic relations for the equation of states are:

$$\begin{aligned} e_f &= e_f(T_f, \upsilon_f) \\ P_f &= P_f(T_f, \upsilon_f) \\ \phi_f &= \upsilon_f \bar{\rho}_f \\ 0 &= n - \sum_{f=1}^{N_f} \upsilon_f \bar{\rho}_f \end{aligned} \tag{10}$$

*Constitutive soil model*

As a result of the explicit MPM formulation, we can derive the constitutive law in the updated Lagrangian framework of "small strain - large deformation". Therefore, the rotation of the particles (representative element volume) is manipulated by rotating the Cauchy stress tensor. First, the deformation gradient is decomposed into the polar rotation tensor $\boldsymbol{R}_s^{n+1}$ and sketch tensor $\boldsymbol{V}_s^{n+1}$ as

$$\boldsymbol{F}_s^{n+1} = \boldsymbol{V}_s^{n+1} \boldsymbol{R}_s^{n+1} \tag{11}$$

Then, before calling the constitutive model, the stress and strain rate tensor are rotated to the reference configuration as

$$\boldsymbol{\sigma}'^{,n*} = (\boldsymbol{R}_s^{n+1})^T \boldsymbol{\sigma}'^{,n*} \boldsymbol{R}_s^{n+1} \tag{12}$$

$$\delta \boldsymbol{\epsilon}^{n*} = (\boldsymbol{R}_s^{n+1})^T \delta \boldsymbol{\epsilon}_s^{n*} \boldsymbol{R}_s^{n+1} \tag{13}$$

Using the constitutive model with the input tensors $\boldsymbol{\sigma}'^{,n*}, \delta\boldsymbol{\epsilon}^{n*}$ to compute the Cauchy stress tensor at the advanced time step $\boldsymbol{\sigma}'^{,n+1*}$ then rotating it back to current configuration

$$\boldsymbol{\sigma}'^{,n+1} = \boldsymbol{R}_s^{n+1} \boldsymbol{\sigma}'^{,n+1*} (\boldsymbol{R}_s^{n+1})^T \tag{14}$$

In this paper, we adopt the hyper-elastic Neo Hooken model for the structure materials and additionally Mohr-Coulomb failure criteria for the soil (porous media) materials. The Cauchy stress of the hyper-elastic Neo Hookean model can be written as:

$$\boldsymbol{\sigma}' = \frac{\lambda ln(J)}{J} + \frac{\mu}{J}(\boldsymbol{F}\boldsymbol{F}^T - \boldsymbol{J}) \tag{15}$$

where $\lambda$ and $\mu$ are bulk and shear modulus ad J is the determinant of the deformation gradient $\boldsymbol{F}$. And the yield function f and flow potentials g of the Mohr-Coulomb can be written as:

$$\begin{aligned} f &= \sigma_1' - \sigma_3' - 2c'cos(\phi') - (\sigma_1' + \sigma_3')sin(\phi') \\ g &= \sigma_1' - \sigma_3' - 2c'cos(\psi') - (\sigma_1' + \sigma_3')sin(\psi') \end{aligned} \tag{16}$$

where the $c'$, $\phi'$ and $\psi'$ are cohesion and friction angle and dilation angle. $\sigma_1'$ and $\sigma_3'$ are maximum and minimum principal stress.

*Turbulent model*

The turbulent effect is modelled using a statistical approach namely large-eddy simulation. In this approach, the micro-scale turbulent influence in the dynamics of the macro-scale motion is computed through simple models like Smagorinsky model. In the Smagorinsky mode, the residual stress tensor is:

$$\tau_{ij} = 2\mu_{eff}(\overline{S}_{ij} - \frac{1}{3}\delta_{ij}\overline{S}_{kk}) + \frac{1}{3}\delta_{ij}\tau_{kk} \tag{17}$$

where the the strain rate tensor is given by

$$\overline{S}_{ij} = \frac{1}{2}(\frac{\delta\overline{\boldsymbol{U}}_i}{\delta x_j} + \frac{\delta\overline{\boldsymbol{U}}_j}{\delta x_i}) \tag{18}$$

and the effective viscosity is sum of molecular viscosity and turbulent viscosity $\mu_{eff} = \mu + \mu_t$ in which the turbulent viscosity $\mu_t$ is calculated by

$$\mu_t = (C_s\triangle)^2\sqrt{2\overline{S}_{ij}\overline{S}_{ij}} \tag{19}$$

where $C_s$ is the Smagorinsky constant and $\triangle = \sqrt[3]{dxdydz}$ is the grid size that defines the subgrid length scale.

*Frictional force for soil-structure interaction*

MPMICE includes a contactlaw for the interaction between soil and structure using the first Coulomb friction contact for MPM presented by Bardenhagen et al. ([19]). The magnitude of the friction force at the contact depends on the friction coefficient $\mu_{fric}$ and the normal force $\boldsymbol{f}_{norm}$ computed from the projection of the contact force in the normal direction.

$$\boldsymbol{f}_{fric} = \mu_{fric}\boldsymbol{f}_{norm} \tag{20}$$

The contact determines whether the soil is sliding or sticking to the structure by comparing the friction force with the sticking force $\boldsymbol{f}_{stick}$ computed from the projection of the contact force in the tangent direction:

$$\begin{aligned} &\text{if } \boldsymbol{f}_{fric} \geq \boldsymbol{f}_{stick} \text{ no sliding} \\ &\text{if } \boldsymbol{f}_{fric} < \boldsymbol{f}_{stick} \text{ sliding occurs} \end{aligned} \tag{21}$$

Frictional sliding between solid materials also generates dissipation and the work rate generated from the sliding can be calculated as:

$$\Delta W_{friction} = \boldsymbol{f}_{fric} d \tag{22}$$

where d is the sliding distance which can computed based on the sliding velocity between two materials.

*Momentum and Energy exchange model*

Currently, the energy exchange coefficient $H_{sf}$ is assumed to be constant for the sake of simplicity. Then the energy exchange can be written as:

$$q_{sf} = H_{sf}(T_f - T_s) \tag{23}$$

On the other hand, the drag force can be calculated as:

$$f_d = K(\boldsymbol{U}_s - \boldsymbol{U}_f) \tag{24}$$

For the momentum exchange between fluid flows and porous media, we assume that the drag force $\boldsymbol{f}_d$ depends on the average grain size of the grains $D_p$, the porosity $n$, the fluid vicosity $\mu_f$, and is propotional to the relative velocities of soil grains and fluid $(\boldsymbol{U}_s - \boldsymbol{U}_f)$. Based on recent investigation of CFD simulations of fluid flow around mono- and bi-disperse packing of spheres for $0.1 < \phi_s < 0.6$ and $Re < 1000$ [20]. The drag force is given by:

$$\boldsymbol{f}_d = \frac{18\phi_s(1-\phi_s)\mu_f}{D_p^2} F(\phi_s, Re)(\boldsymbol{U}_s - \boldsymbol{U}_f) \tag{25}$$

where Reynolds number $Re$ are computed as:

$$Re = \frac{n\rho_f D_p}{\mu_f} \| \|(\boldsymbol{U}_s - \boldsymbol{U}_f)\| \tag{26}$$

The function $F(\phi_s, Re)$ can be calculated as:

$$F(\phi_s, Re) = F(\phi_s, 0) + \frac{0.413Re}{24(1-\phi_s)^2}\left(\frac{(1-\phi_s)^{-1} + 3\phi_s(1-\phi_s) + 8.4Re^{-0.343}}{1 + 10^{3\phi_s}Re^{-(1+4\phi_s)/2}}\right) \tag{27}$$

where the low Reynold coefficient $F(\phi_s, Re \rightarrow 0)$ is:

$$F(\phi_s, 0) = \frac{10\phi_s}{(1-\phi_s)^2} + (1-\phi_s)^2(1 + 1.5\sqrt{\phi_s}) \tag{28}$$

When validating the model with analytical solution, it requires to know the hydraulic conductivity. In such case, we convert the equation (25) to Kozeny-Carman formula by assuming $F(\phi_s, Re) = 10\phi_s/(1-\phi_s)^2$, then the hydraulic conductivity will be expressed as $D_p^2(1-\phi_s)^3/180\mu\phi_s^2$.

*Solving momentum and energy exchange with an implicit solver*

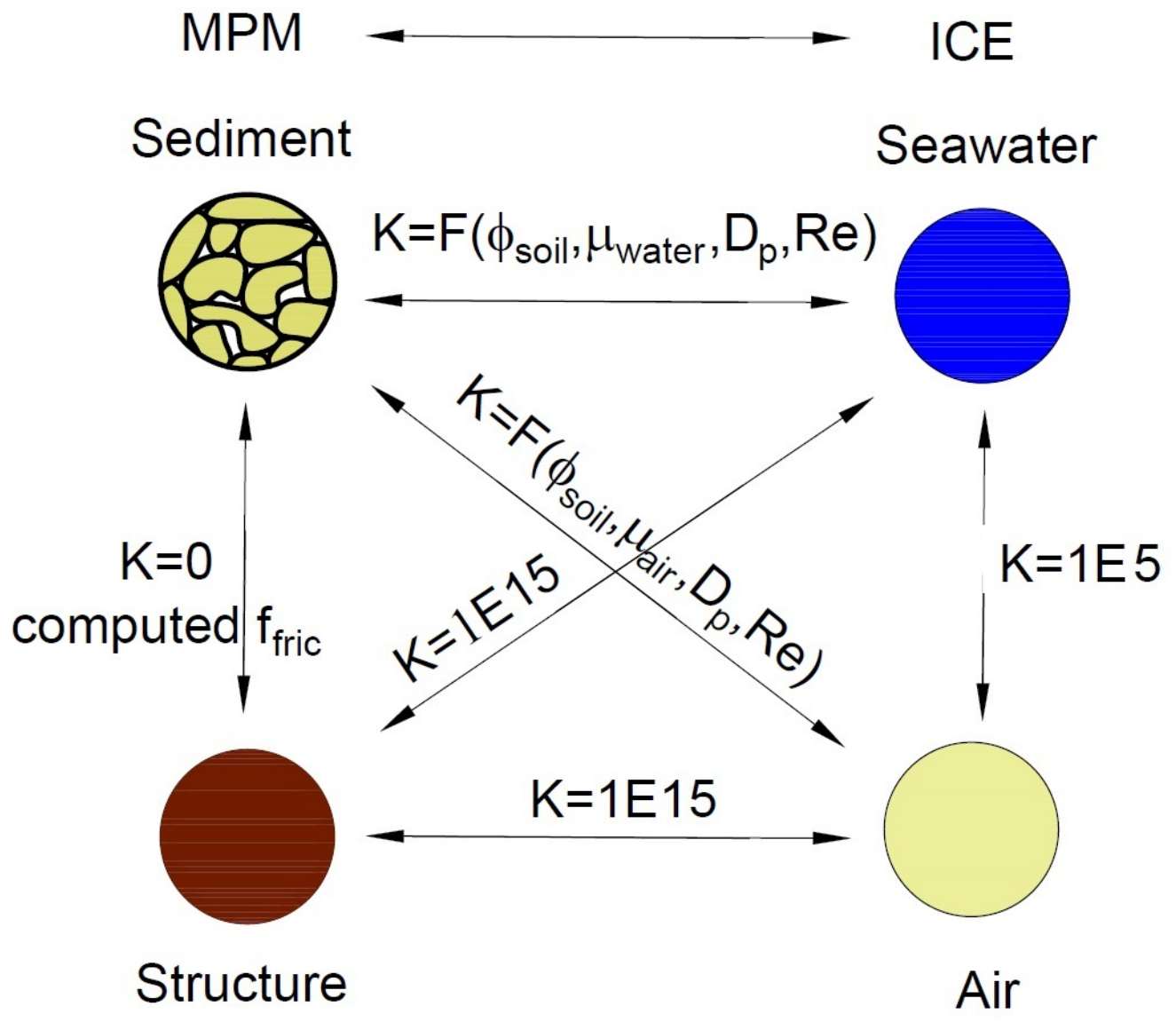

*Figure 3: Momentum exchange coefficient between materials*

The derivation of the implicit integration for the momentum exchange is presented in the Appendix's section 'Momentum and energy exchange with an implicit solver'. The linear equations for multi phases i,j=1:$N$ has the

form:

$$\begin{vmatrix} (1+\beta_{ij}) & -\beta_{ij} \\ -\beta_{ji} & (1+\beta_{ji}) \end{vmatrix} \begin{vmatrix} \Delta \boldsymbol{U}_i \\ \Delta \boldsymbol{U}_j \end{vmatrix} = \begin{vmatrix} \beta_{ij}(\boldsymbol{U}_i^* - \boldsymbol{U}_j^*) \\ \beta_{ji}(\boldsymbol{U}_i^* - \boldsymbol{U}_j^*) \end{vmatrix}$$

where the intermediate velocity for fluid phases f=1:$N_f$ and for solid/porous phases s=1:$N_s$ can be calculated by

$$\begin{aligned} \boldsymbol{U}_f^* &= \boldsymbol{U}_f^n + \Delta t(-\frac{\nabla P_f^{n+1}}{\rho_f^n} + \frac{\nabla \cdot \boldsymbol{\tau}_f^n}{\overline{\rho}_f^n} + \boldsymbol{b}) \\ \boldsymbol{U}_s^* &= \boldsymbol{U}_s^n + \Delta t(\frac{\nabla \cdot \boldsymbol{\sigma}'^n}{\overline{\rho}_s^n} - \frac{\nabla P_f^{n+1}}{\rho_s} + \boldsymbol{b}) \end{aligned} \tag{29}$$

Also, the momentum exchange coefficient can be computed at every time step as $\beta_{12} = K/\overline{\rho}_f^n$ and $\beta_{21} = K/\overline{\rho}_s^n$ with the coefficient depending on the different type of interactions (see Figure 3) as for example:

1. The drag force is set to zero in soil-structure interactions, and instead the frictional force is computed.
2. As a result of fluid-structure interaction, the momentum exchange coefficient should be extremely high when the solid material points are considered to be zero-porosity/zero-permeability.
3. In the case of soil-fluid interaction, the drag force is calculated using the equation (25). Considering that air has a much lower viscosity than water, its drag force is much lower than the drag force of water in a pore.
4. A momentum exchange coefficient of 1E5 is applied between multiphase flows. This value is far higher than reality [21], but it is necessary to have enough numerical stability to conduct simulations in the numerical example.

Similar approach applied for the ernergy exchange term leading to:

$$\begin{vmatrix} (1+\eta_{ij}) & -\eta_{ij} \\ -\eta_{ji} & (1+\eta_{ji}) \end{vmatrix} \begin{vmatrix} \Delta T_i \\ \Delta T_j \end{vmatrix} = \begin{vmatrix} \eta_{ij}(T_i^n - T_j^n) \\ \eta_{ji}(T_j^n - T_i^n) \end{vmatrix}$$

with $\eta$ is the energy exchange coefficient.

*Equation of state for fluid plases*

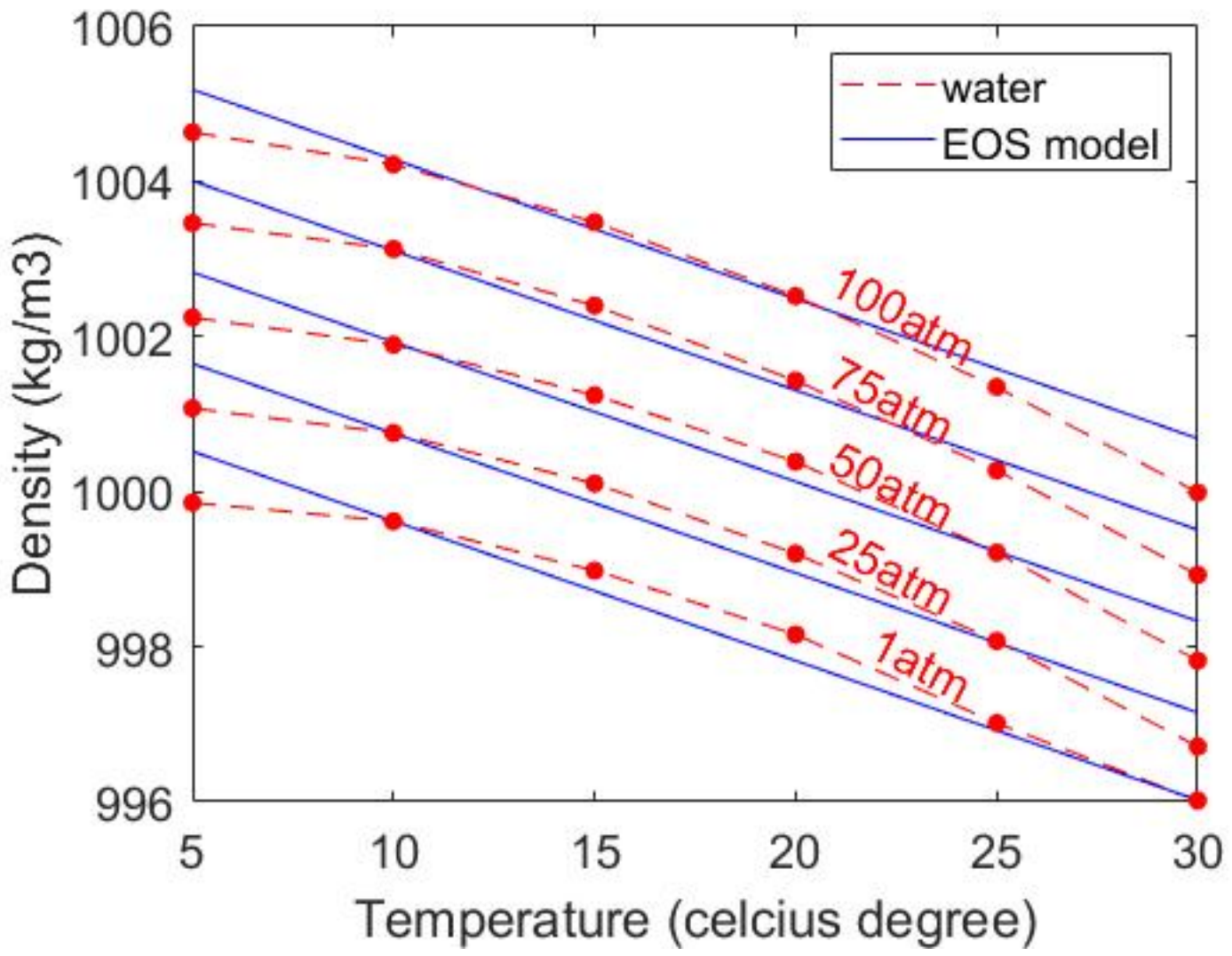

*Figure 4: Equation of state of water*

The equation of state establishes relations between thermodynamics variables $[P_f, \rho_f, T_f]$. The choice of the equation of state depends on the types of the fluid materials. For example, for the air, it is possible to assume the equation of state for the perfect gas which obeys:

$$P_f = \rho_f R T_f \tag{30}$$

where $R$ is the gas constant. For the water, a simple linear equation of state is in the following form:

$$P_f = P_{ref} + K_f(\rho_f - \rho_{ref} - \alpha_f(T_f - T_{ref})) \tag{31}$$

where reference pressure $P_{ref} = 1$ atm $= 101325$ Pa, reference temperature $T_{ref} = 10°\text{C}$, reference density $\rho_{ref} = 999.8$ kg/m3, the bulk modulus of water $K_f = 2$ GPa, and the water thermal expansion $\alpha_f = 0.18$ $°C^{-1}$. Equation (31) matches well with the state of the water (see Figure 4).

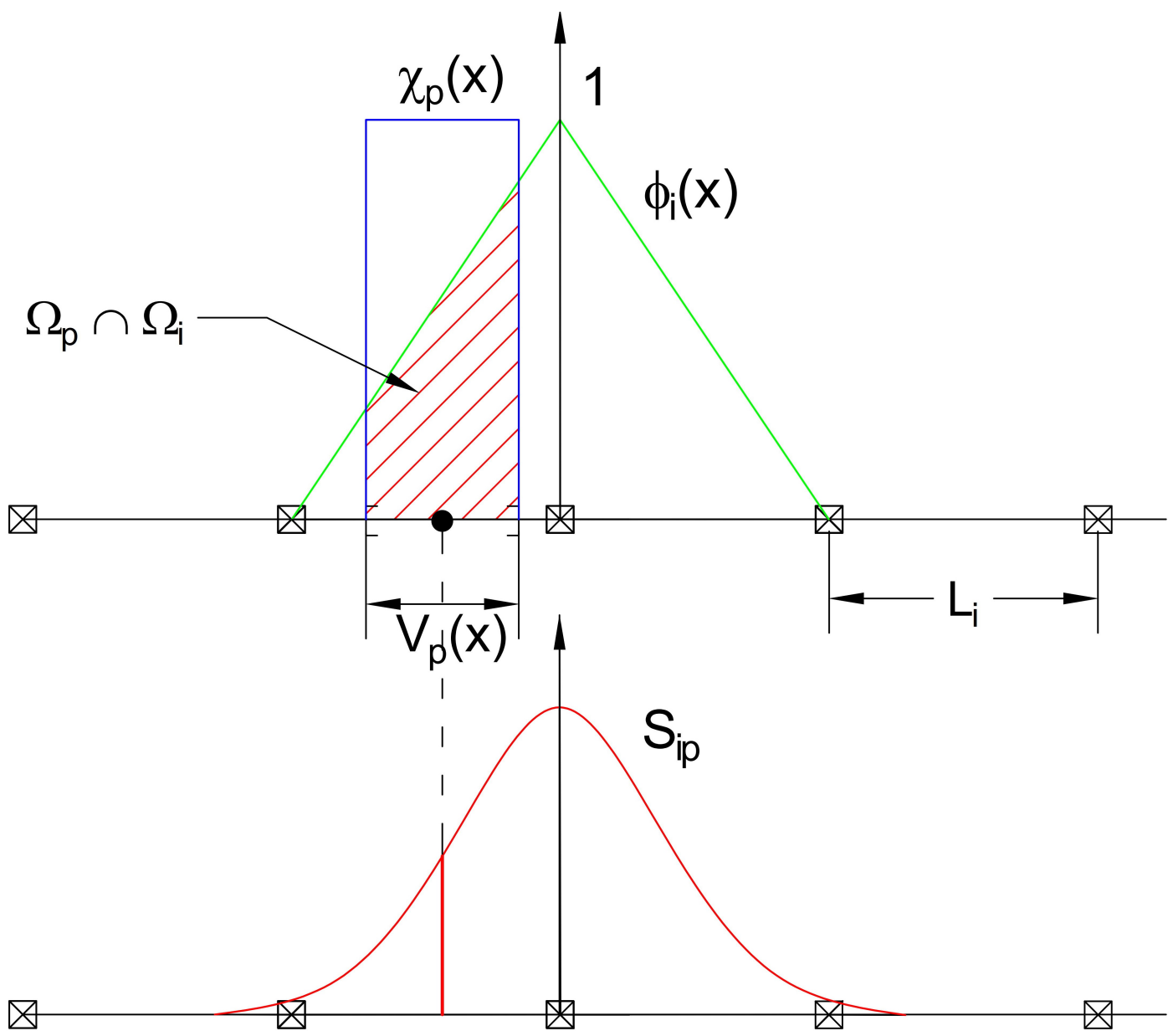

*Figure 5: GIMP weighting function (red as a convolution of the linear basis shape function (green) and the charateristic function (blue)*

The fluid phases are discretized in the grid with the state variables stored at the centroid of the cells $[\rho_{fc}, \boldsymbol{U}_{fc}, T_{fc}, \upsilon_{fc}]$ while the solid phase is discretized in the particles with the state variables $[m_{sp}, \boldsymbol{U}_{sp}, T_{sp}, \boldsymbol{\sigma}'_{sp}]$. In the Material Point Method, we use the generalized interpolation technique [22] using the weight function as a convolution of a grid shape function $N_i(\boldsymbol{x})$ in a nodal domain $\Omega_i$ and a characteristic function $\chi_p(\boldsymbol{x})$ in a particle domain $\Omega_p$ with the volume $V_p(\boldsymbol{x})$ as follows:

$$S_{ip} = \frac{1}{V_p} \int_{\Omega_i \cap \Omega_p} N_i(\boldsymbol{x}) \chi_p(\boldsymbol{x}) d\boldsymbol{x} \tag{32}$$

where the volume $V_p(\boldsymbol{x})$ of the material point p can be calculated as:

$$V_p = \int_{\Omega_p} \chi_p(\boldsymbol{x}) d\boldsymbol{x} \tag{33}$$

The charateristic function is the Heaviside function as $\chi_p = 1$ if $\boldsymbol{x} \in \Omega_p$, otherwise 0 (see Figure 5). For the interpolation of the centroid of the cell, the linear basis function is used as:

$$S_{ci} = N_i(\boldsymbol{x_c}) \tag{34}$$

The time discretization are solved using the following steps:

*Interpolation from Solid Particle to Grid*

The nodal values of the solid state (mass, velocity, temperature, volume) are:

$$\begin{aligned}
m^n_{si} &= \sum S_{ip} m_{sp} \\
\boldsymbol{U}^n_{si} &= \frac{\sum S_{ip}(m\boldsymbol{U})^n_{sp}}{m^n_{si}} \\
T^n_{si} &= \frac{\sum S_{ip}(mT)^n_{sp}}{m^n_{si}} \\
V^n_{si} &= \frac{\sum S_{ip}(mV)^n_{sp}}{m^n_{si}} \\
\boldsymbol{\sigma}^n_{si} &= \frac{\sum S_{ip}(\boldsymbol{\sigma}V)^n_{sp}}{V^n_{si}}
\end{aligned} \tag{35}$$

The nodal internal forces is calculated by

$$\boldsymbol{f}^{int,n}_{si} = -\sum \nabla S_{ip}(\sigma'_{sp})^n V^n_{sp} \tag{36}$$

The nodal external forces $f^{ext,n}_{si}$ and extra momentum from contact forces are computed here. The nodal velocity and nodal temperature are applied boundary conditions.
Then we compute the solid cell variables as:

$$\begin{aligned}
m^n_{sc} &= \sum S_{ci} m_{si} \\
\rho^n_{sc} &= \frac{m^n_{sc}}{V} \\
\boldsymbol{U}^n_{sc} &= \sum S_{ci}\boldsymbol{U}^n_{si} \\
T^n_{sc} &= \sum S_{ci} T^n_{si} \\
V^n_{sc} &= \sum S_{ci} V^n_{si} \\
\boldsymbol{\sigma}^n_{sc} &= \sum S_{ci}\boldsymbol{\sigma}^n_{si}
\end{aligned} \tag{37}$$

*Compute equation of state for fluid phase*

Considering the total fluid materal volume of a cell is:

$$V_{total} = \sum_{f=1}^{N_f} M_f \upsilon_f \tag{38}$$

Then we need to find $P_{eq}$ which allows each fluid materials obey their equation of states $[P_f, \rho_f, \upsilon_f, T_f, e_f]$ but also allow mass of all fluid materials to fill the entire the pore volume without ongoing compression or expansion following the condition:

$$0 = n - \sum_{f=1}^{N_f} \upsilon_f \overline{\rho}_f \tag{39}$$

Then, we can use he Newton-Raphson interation to find the value of $P_{eq}$ which satisfies the equation (38, 39) and each equation of states of each fluid materials.

*Compute faced-centered velocity*

Following the derivation in the Appendix: Advanced Fluid Pressure, we first compute the fluid face-centered velocity as

$$\boldsymbol{U}^*_{f,FC} = \frac{(\overline{\rho}\boldsymbol{U})^n_{f,FC}}{\overline{\rho}^n_{f,FC}} + \Delta t(-\frac{\nabla^{FC} P_{eq}}{\rho^n_{f,FC}} + \frac{\nabla^{FC} \cdot \boldsymbol{\tau}^n}{\overline{\rho}_{s,FC}} + \boldsymbol{b}) \tag{40}$$

The equation (40) is discretized in three dimension (noted that $\nabla^{FC} \cdot \boldsymbol{\tau} = 0$), for example the discretized equation in the x direction is

$$U^*_{fx} = \frac{(\overline{\rho}U)^n_{fx,R} + (\overline{\rho}U)^n_{fx,L}}{\overline{\rho}^n_{fx,L} + \overline{\rho}^n_{fx,R}} + \Delta t(-\frac{2(\upsilon^n_{fx,L}\upsilon^n_{fx,R})}{\upsilon^n_{fx,L} + \upsilon^n_{fx,R}} \frac{P_{eqx,R} - P_{eqx,L}}{\Delta x} + b_x) \tag{41}$$

The face-centered solid velocity can be calculated as

$$\boldsymbol{U}^*_{s,FC} = \frac{(\overline{\rho}\boldsymbol{U})^n_{s,FC}}{\overline{\rho}^n_{s,FC}} + \Delta t(\frac{\nabla^{FC} \cdot \boldsymbol{\sigma}'^n_c}{\overline{\rho}_{s,FC}} - \frac{\nabla^{FC} P_{eq}}{\rho_s} + \boldsymbol{b}) \tag{42}$$

The equation (42) is discretized in three dimension(noted that $\nabla^{FC} \cdot \sigma_{ij} = 0$ with $i \neq j$), for example the discretized equation in the x direction is

$$U^*_{sx} = \frac{(\overline{\rho}U)^n_{sx,R} + (\overline{\rho}U)^n_{sx,L}}{\overline{\rho}^n_{sx,L} + \overline{\rho}^n_{sx,R}} + \Delta t(\frac{2(\sigma_{xx,R} - \sigma_{xx,L})}{(\overline{\rho}^n_{sx,L} + \overline{\rho}^n_{sx,R})\Delta x} - \frac{P_{eqx,R} - P_{eqx,L}}{\rho_s \Delta x} + b_x) \tag{43}$$

Computing the modified faced-centered velocity $\boldsymbol{U}^{L}_{FC}$ considering the momentum exchange (see the Appendix: Momentum exchange with an implicit solve)

$$\begin{aligned} \boldsymbol{U}^{L}_{f,FC} &= \boldsymbol{U}^{*}_{f,FC} + \Delta\boldsymbol{U}_{f,FC} \\ \boldsymbol{U}^{L}_{s,FC} &= \boldsymbol{U}^{*}_{s,FC} + \Delta\boldsymbol{U}_{s,FC} \end{aligned} \tag{44}$$

Solving the linear equation below to obtain the increment of velocity with i,j $= 1 : N$ as:

$$\begin{vmatrix} (1+\beta_{ij}) & -\beta_{ij} \\ -\beta_{ji} & (1+\beta_{ji}) \end{vmatrix} \begin{vmatrix} \Delta\boldsymbol{U}_{i,FC} \\ \Delta\boldsymbol{U}_{j,FC} \end{vmatrix} = \begin{vmatrix} \beta_{ij}(\boldsymbol{U}^{*}_{i,FC} - \boldsymbol{U}^{*}_{j,FC}) \\ \beta_{ji}(\boldsymbol{U}^{*}_{j,FC} - \boldsymbol{U}^{*}_{i,FC}) \end{vmatrix}$$

*Compute faced-centered temperature*

Similar to the velocity, the faced temperature is computed as:

$$T^{n}_{fx} = \frac{(\overline{\rho}T)^{n}_{fx,R} + (\overline{\rho}T)^{n}_{fx,L}}{\overline{\rho}^{n}_{fx,L} + \overline{\rho}^{n}_{fx,R}} \tag{45}$$

Computing the modified faced-centered temperature $T^{L}_{FC}$ considering the energy exchange (see the Appendix: Momentum and energy exchange with an implicit solver)

$$\begin{aligned} T^{L}_{f,FC} &= T^{n}_{f,FC} + \Delta T_{f,FC} \\ T^{L}_{s,FC} &= T^{n}_{s,FC} + \Delta T_{s,FC} \end{aligned} \tag{46}$$

Solving the linear equation below to obtain the increment of velocity with i,j $= 1 : N$ as:

$$\begin{vmatrix} (1+\eta_{ij}) & -\eta_{ij} \\ -\eta_{ji} & (1+\eta_{ji}) \end{vmatrix} \begin{vmatrix} \Delta T_{i,FC} \\ \Delta T_{j,FC} \end{vmatrix} = \begin{vmatrix} \eta_{ij}(T^{n}_{i,FC} - T^{n}_{j,FC}) \\ \eta_{ji}(T^{n}_{j,FC} - T^{n}_{i,FC}) \end{vmatrix}$$

*Compute fluid pressure (implicit scheme)*

For single phase flow, the increment of the fluid pressure can be computed as:

$$\kappa_f \frac{\Delta P}{dt} = \nabla^{c} \cdot \boldsymbol{U}^{n+1}_{f,FC} \tag{47}$$

For multi-phase flows, the increment of the fluid pressure of the mixture can be computed as:

$$\kappa \frac{\Delta P}{dt} = \sum_{f=1}^{N_f} \nabla^{c} \cdot (\phi_{f,FC}\boldsymbol{U}_{f,FC})^{n+1} \tag{48}$$

where $\kappa = \sum_{f=1}^{N_f}(\phi_{f,FC}\kappa_f)/\sum_{f=1}^{N_f}(\phi_{f,FC})$. Then, the fluid pressure at cell center is:

$$P_c^{n+1} = P_{eq} + \Delta P_c^n \tag{49}$$

Finally, the faced-centered advanced fluid pressure is

$$P_{FC}^{n+1} = (\frac{P_{c,L}^{n+1}}{\overline{\rho}_{f,L}^n} + \frac{P_{c,R}^{n+1}}{\overline{\rho}_{f,R}^n})/(\frac{1}{\overline{\rho}_{f,L}^n} + \frac{1}{\overline{\rho}_{f,R}^n}) = (\frac{P_{c,L}^{n+1}\overline{\rho}_{f,R}^n + P_{c,R}^{n+1}\overline{\rho}_{f,L}^n}{\overline{\rho}_{f,L}^n\overline{\rho}_{f,R}^n}) \tag{50}$$

*Compute viscous shear stress term of the fluid phase*

This part compute the viscous shear stress $\Delta(m\boldsymbol{U})_{fc,\tau}$ for a single vicous compressible Newtonian fluid and optionally shear stress induced by the turbulent model.

*Compute nodal internal temperature of the solid phase*

The nodal internal temperature rate is computed based on the heat conduction model

$$dT_{si}^L = \frac{(\Delta W_{si}^n + \Delta W_{fric,i}^n + \nabla^i \cdot \boldsymbol{q}_{si}^n)}{m_{si}^n c_v} \tag{51}$$

where $\Delta W_{si}^n = \boldsymbol{\sigma'} : \frac{D_s(\boldsymbol{\epsilon}_s^p)}{Dt}$ is the mechanical work rate computed from the constitutive model with $\boldsymbol{\epsilon}_s^p$ is the plastic strain, $\Delta W_{fric,i}^n$ is the work rate compted from the contact law due to the frictional sliding between solid materials. The heat flux is $\boldsymbol{q}_s = \overline{\rho}_s\beta_s\nabla T_s$ with $\beta_s$ being the thermal conductivity of the solid materials.

$$T_{si}^L = T_{si}^n + dT_{si}^L \tag{52}$$

*Compute and integrate acceleration of the solid phase*

After interpolating from material points to the nodes, the nodal acceleration and velocity are calculate by

$$\boldsymbol{a}_{si}^{L-} = \frac{\boldsymbol{f}_{si}^{int,n} + \boldsymbol{f}_{si}^{ext,n}}{m_{si}^n} + \boldsymbol{g} \tag{53}$$

$$\boldsymbol{U}_{si}^{L-} = \boldsymbol{U}_{si}^n + \boldsymbol{a}_{si}^{L-}\Delta t \tag{54}$$

*Compute Lagrangian value (mass, momentum and energy)*

For the fluid phase, the linear momentum rate, the energy rate are

$$\Delta(m\boldsymbol{U})_{fc} = Vn_c^n \nabla^c P_c^{n+1} + \Delta(m\boldsymbol{U})_{fc,\tau} + V\overline{\rho}_{fc}^n g \tag{55}$$

$$\Delta(me)_{fc} = Vn_c^n P_c^{n+1} \nabla^c \cdot \boldsymbol{U}_{f,FC}^* + \nabla^c \cdot \boldsymbol{q}_{fc}^n \tag{56}$$

The Lagrangian value of the mass, linear momentum and energy of fluid phases without momentum exchange are

$$m_{fc}^L = V\overline{\rho}_{fc}^n \tag{57}$$

$$(m\boldsymbol{U})_{fc}^{L-} = V\overline{\rho}_{fc}^n \boldsymbol{U}_{fc}^n + \Delta(m\boldsymbol{U})_{fc} \tag{58}$$

$$(me)_{fc}^{L-} = V\overline{\rho}_{fc}^n T_{fc}^n c_v + \Delta(me)_{fc} \tag{59}$$

For the solid phase, the Lagrangian value of the linear momentum and energy of solid phase are

$$m_{sc}^L = m_{sc}^n \tag{60}$$

$$(m\boldsymbol{U})_{sc}^{L-} = \sum S_{ci} m_{si}^n \boldsymbol{U}_{si}^{L-} + V(1-n_c^n)\nabla^c P_{fc}^{n+1} \tag{61}$$

$$(me)_{sc}^{L-} = \sum S_{ci} m_{si}^n T_{si}^L c_v \tag{62}$$

To consider the momentum exchange, the Lagrangian velocity is modified as

$$\begin{aligned} \boldsymbol{U}_{fc}^L &= \boldsymbol{U}_{fc}^{L-} + \Delta\boldsymbol{U}_{fc} \\ \boldsymbol{U}_{sc}^L &= \boldsymbol{U}_{sc}^{L-} + \Delta\boldsymbol{U}_{sc} \end{aligned} \tag{63}$$

where the cell-centered intermediate velocity can be calculated by

$$\begin{aligned} \boldsymbol{U}_{fc}^{L-} &= \frac{(m\boldsymbol{U})_{fc}^{L-}}{m_{fc}^L} \\ \boldsymbol{U}_{sc}^{L-} &= \frac{(m\boldsymbol{U})_{sc}^{L-}}{m_{sc}^L} \end{aligned} \tag{64}$$

And the increment of the velocity $\boldsymbol{U}_{fc}$, $\Delta\boldsymbol{U}_{sc}$ can be computed by solving the linear equation with i,j = 1 : $N$ as:

$$\begin{vmatrix} (1+\beta_{ij}) & -\beta_{ij} \\ -\beta_{ji} & (1+\beta_{ji}) \end{vmatrix} \begin{vmatrix} \Delta\boldsymbol{U}_{i,c} \\ \Delta\boldsymbol{U}_{j,c} \end{vmatrix} = \begin{vmatrix} \beta_{ij}(\boldsymbol{U}_{i,c}^* - \boldsymbol{U}_{j,c}^*) \\ \beta_{ji}(\boldsymbol{U}_{j,c}^* - \boldsymbol{U}_{i,c}^*) \end{vmatrix}$$

To consider the energy exchange, the Lagrangian temperature is modified as

$$\begin{aligned} T_{fc}^{L} &= T_{fc}^{L-} + \Delta T_{fc} \\ T_{sc}^{L} &= T_{sc}^{L-} + \Delta T_{sc} \end{aligned} \tag{65}$$

where the cell-centered intermediate temperature can be calculated by

$$\begin{aligned} T_{fc}^{L-} &= \frac{(mT)_{fc}^{L-}}{m_{fc}^{L} c_v} \\ T_{sc}^{L-} &= \frac{(mT)_{sc}^{L-}}{m_{sc}^{L} c_v} \end{aligned} \tag{66}$$

And the increment of the velocity can be computed by solving the linear equation with i,j = $1 : N$ as:

$$\begin{vmatrix} (1+\eta_{ij}) & -\eta_{ij} \\ -\eta_{ji} & (1+\eta_{ji}) \end{vmatrix} \begin{vmatrix} \Delta T_{i,c} \\ \Delta T_{j,c} \end{vmatrix} = \begin{vmatrix} \eta_{ij}(T_{i,c}^{n} - T_{j,c}^{n}) \\ \eta_{ji}(T_{j,c}^{n} - T_{i,c}^{n}) \end{vmatrix}$$

Finally, we obtain the cell-centered solid acceleration and temperature rate as

$$d\boldsymbol{U}_{sc}^{L} = \frac{(m\boldsymbol{U})_{sc}^{L} - (m\boldsymbol{U})_{sc}^{n}}{m_{sc}^{L}\Delta t} \tag{67}$$

$$dT_{sc}^{L} = \frac{(me)_{sc}^{L} - (me)_{sc}^{n}}{m_{sc}^{L} c_v \Delta t} \tag{68}$$

*Compute Lagrangian specific volume of the fluid phase*

To compute the Lagrangian value of the specific volume of the fluid phase, we need to compute the Lagrangian temperature rate as below

$$T_{fc}^{n+1} = \frac{(me)_{fc}^{L}}{m_{fc}^{L} c_v} \tag{69}$$

$$\frac{D_f T_{fc}}{Dt} = \frac{T_{fc}^{n+1} - T_{fc}^{n}}{\Delta t} \tag{70}$$

As such, the Lagrangian specific volume rate is:

$$\Delta(m\upsilon)_{fc} = V f_{fc}^{\phi} \nabla \cdot \boldsymbol{U} + (\phi_{fc}\alpha_{fc}\frac{D_f T_{fc}}{Dt} - f_{fc}^{\phi}\sum_{n=1}^{N}\phi_{nc}\alpha_{nc}\frac{D_n T_{nc}}{Dt}) \tag{71}$$

where $f_f^{\phi} = (\phi_f \kappa_f)/(\sum_{n=1}^{N}\phi_n\kappa_n)$ and $\boldsymbol{U} = \nabla \cdot (\sum_{s=1}^{N_s}\phi_{sc}\boldsymbol{U}_{sc} + \sum_{f=1}^{N_f}\phi_{fc}\boldsymbol{U}_{fc})$. Finally, the Lagrangian specific volume is

$$(m\upsilon)_{fc}^{L} = V\overline{\rho}_{f,c}^{n}\upsilon_{fc}^{n} + \Delta(m\upsilon)_{fc} \tag{72}$$

*Compute advection term and advance in time*

The time advanced mass, linear momentum, energy and specific volume are:

$$m_{fc}^{n+1} = m_{fc}^L - \Delta t \nabla \cdot (\overline{\rho}_{fc}^L, \boldsymbol{U}_{f,FC}^L) \tag{73}$$

$$(m\boldsymbol{U})_{fc}^{n+1} = (m\boldsymbol{U})_{fc}^L - \Delta t \nabla \cdot ((\overline{\rho}\boldsymbol{U})_{fc}^L, \boldsymbol{U}_{f,FC}^L) \tag{74}$$

$$(me)_{fc}^{n+1} = (me)_{fc}^L - \Delta t \nabla \cdot ((\overline{\rho} c_v T)_{fc}^L, \boldsymbol{U}_{f,FC}^L) \tag{75}$$

$$(m\upsilon)_{fc}^{n+1} = (m\upsilon)_{fc}^L - \Delta t \nabla \cdot ((\overline{\rho}\upsilon)_{fc}^L, \boldsymbol{U}_{f,FC}^L) \tag{76}$$

Finally, the state variables of the fluid phases of the next time step are

$$\overline{\rho}_{fc}^{n+1} = \frac{m_{fc}^{n+1}}{V} \tag{77}$$

$$\boldsymbol{U}_{fc}^{n+1} = \frac{(m\boldsymbol{U})_{fc}^{n+1}}{m_{fc}^{n+1}} \tag{78}$$

$$T_{fc}^{n+1} = \frac{(me)_{fc}^{n+1}}{m_{fc}^{n+1}} \tag{79}$$

$$\upsilon_{fc}^{n+1} = \frac{(m\upsilon)_{fc}^{n+1}}{m_{fc}^{n+1}} \tag{80}$$

*Interpolate from cell to node of the solid phase*

First we interpolate the acceleration, velocity and temperature rate to the node

$$\boldsymbol{a}_{si}^n = \sum S_{ci} d\boldsymbol{U}_{sc}^L \tag{81}$$

$$\boldsymbol{U}_{si}^{n+1} = \sum S_{ci} d\boldsymbol{U}_{sc}^L \Delta t \tag{82}$$

$$dT_{si}^n = \sum S_{ci} dT_{sc}^L \tag{83}$$

Then the boundary condition and contact forces are applied to the nodal velocity and the acceleration is modified by

$$\boldsymbol{a}_{si}^n = \frac{\boldsymbol{v}_{si}^{n+1} - \boldsymbol{v}_{si}^n}{\Delta t} \tag{84}$$

*Update the particle variables*

The state variables of the solid phase $[\boldsymbol{U}_{sp}^{n+1}, \boldsymbol{x}_{sp}^{n+1}, \nabla\boldsymbol{U}_{sp}^{n+1}, T_{sp}^{n+1}, \nabla T_{sp}^{n+1}, \boldsymbol{F}_{sp}^{n+1}, V_{sp}^{n+1}]$ (velocity, position, velocity gradient, temperature, temperature gradient, deformation gradient, volume) are updated here

$$\boldsymbol{U}_{sp}^{n+1} = \boldsymbol{U}_{sp}^{n} + \sum S_{sp}\boldsymbol{a}_{si}^{n}\Delta t \tag{85}$$

$$\boldsymbol{x}_{sp}^{n+1} = \boldsymbol{x}_{sp}^{n} + \sum S_{sp}\boldsymbol{U}_{si}^{n+1}\Delta t \tag{86}$$

$$\nabla\boldsymbol{U}_{sp}^{n+1} = \sum \nabla S_{sp}\boldsymbol{U}_{si}^{n+1} \tag{87}$$

$$T_{sp}^{n+1} = T_{sp}^{n} + \sum S_{sp}dT_{si}^{n}\Delta t \tag{88}$$

$$\nabla T_{sp}^{n+1} = \sum \nabla S_{sp}T_{sp}^{n}\Delta t \tag{89}$$

$$\boldsymbol{F}_{sp}^{n+1} = (\boldsymbol{I} + \nabla\boldsymbol{U}_{sp}^{n+1}\Delta t)\boldsymbol{F}_{sp}^{n} \tag{90}$$

$$V_{sp}^{n+1} = det(\boldsymbol{F}_{sp}^{n+1})V_{sp}^{o} \tag{91}$$

Finally, the effective stress $(\boldsymbol{\sigma'})^{n+1}$ is updated from the constitutive model and the pore water pressure is interpolated from the cell as:

$$p_f^{n+1} = \sum S_{si}P_c^{n+1} \tag{92}$$

## Numerical examples

All input files and the analytical calculations in this section are provided in the Github repository for the reproduction of the numerical results.

*Fluid Flow through isothermal porous media*

Fluid flow through porous media is important in many engineering disciplines, like predicting water flow in soil. Fluid flow velocity in one dimension can be calculated from the porous media's hydraulic conductivity $K$ as:

$$U_f = K\frac{\Delta p_f}{L} \tag{93}$$

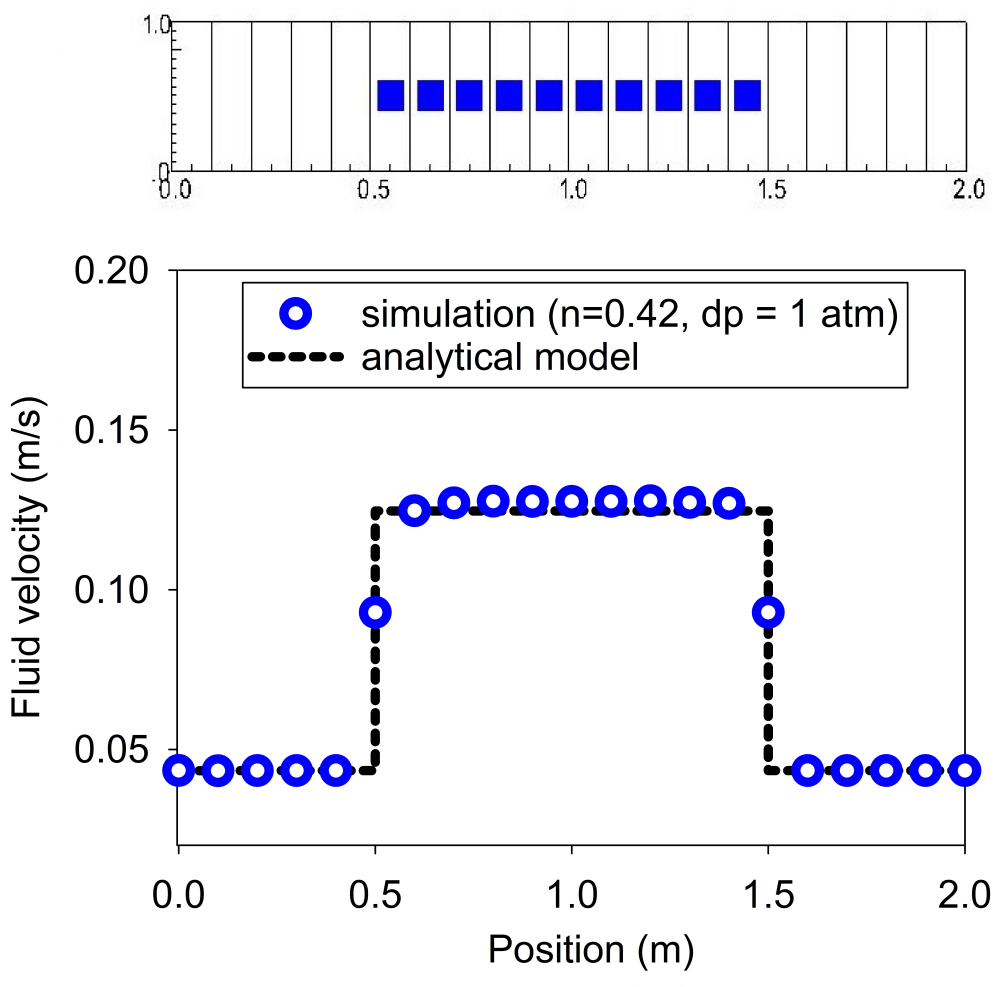

*(a) Discretization of the model*

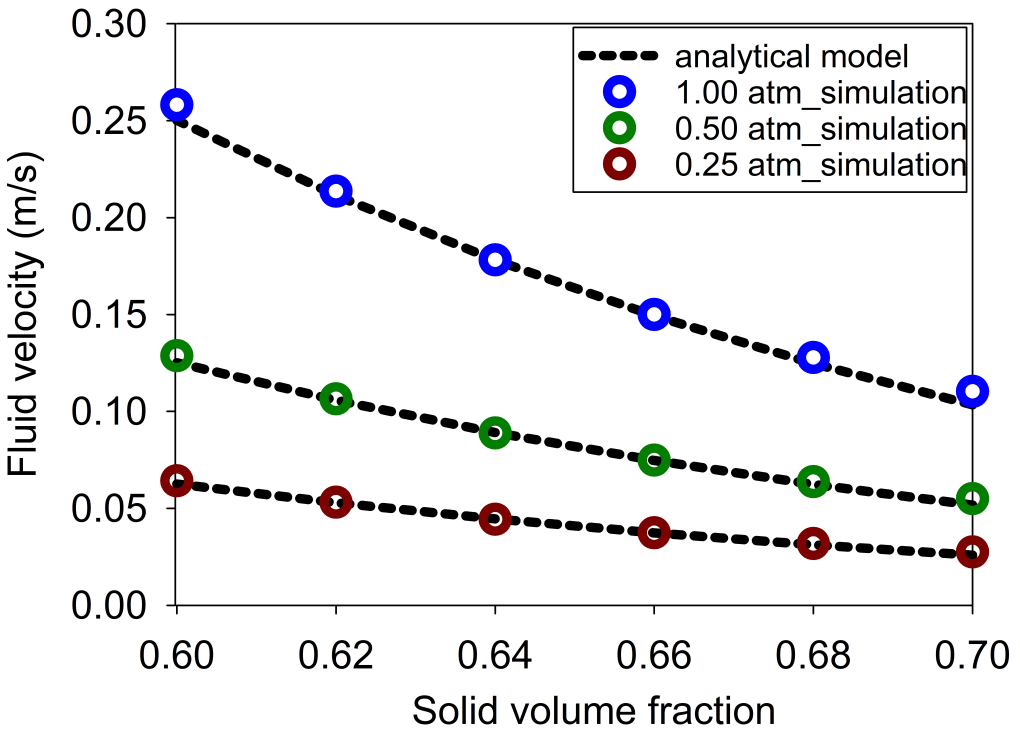

*(b) Comparision between analytical model and simulation*

*Figure 6: Numerical results of the fluid flow through isothermal porous media*

If the Carman-Kozeny formula is adopted $F = 10\phi_s/(1-\phi_s)^2$, the hydraulic conductivity will be expressed as $K = d^2(1-\phi_s)^3/180\mu\phi_s^2$. Then, the analytical formula of average velocity in one dimension through the porous media is:

$$U_f = \frac{1}{n}\frac{d^2(1-\phi_s)^3}{180\mu\phi_s^2}\frac{\Delta p_f}{L} \tag{94}$$

Our numerical model is validated by modeling fluid flow through a 1m long porous media. This fluid has water properties (bulk modulus is 2GPa, density is 998 kg/m3 at 5 degrees Celsius and 10325 Pa (1atm) pressure, dynamic viscosity $\mu$ is 1mPa s). The porous media is modeled by elastic material with Young's modulus is 10 MPa, Poisson's ratio is 0.3, and density is 2650 kg/m3. The volume fraction of porous media $\phi_s$ is [0.6, 0.62, 0.66, 0.68, 0.7] and the average grain diameter $d$ is 1mm. The model is discretized in 20 finite element and the porous media in 10 finite element with 1 material point per element. The pressure gradient is applied with three different value [0.25, 0.5, 1] atm. Figure 6 shows a good agreement of fluid flow prediction between the theory and the model.

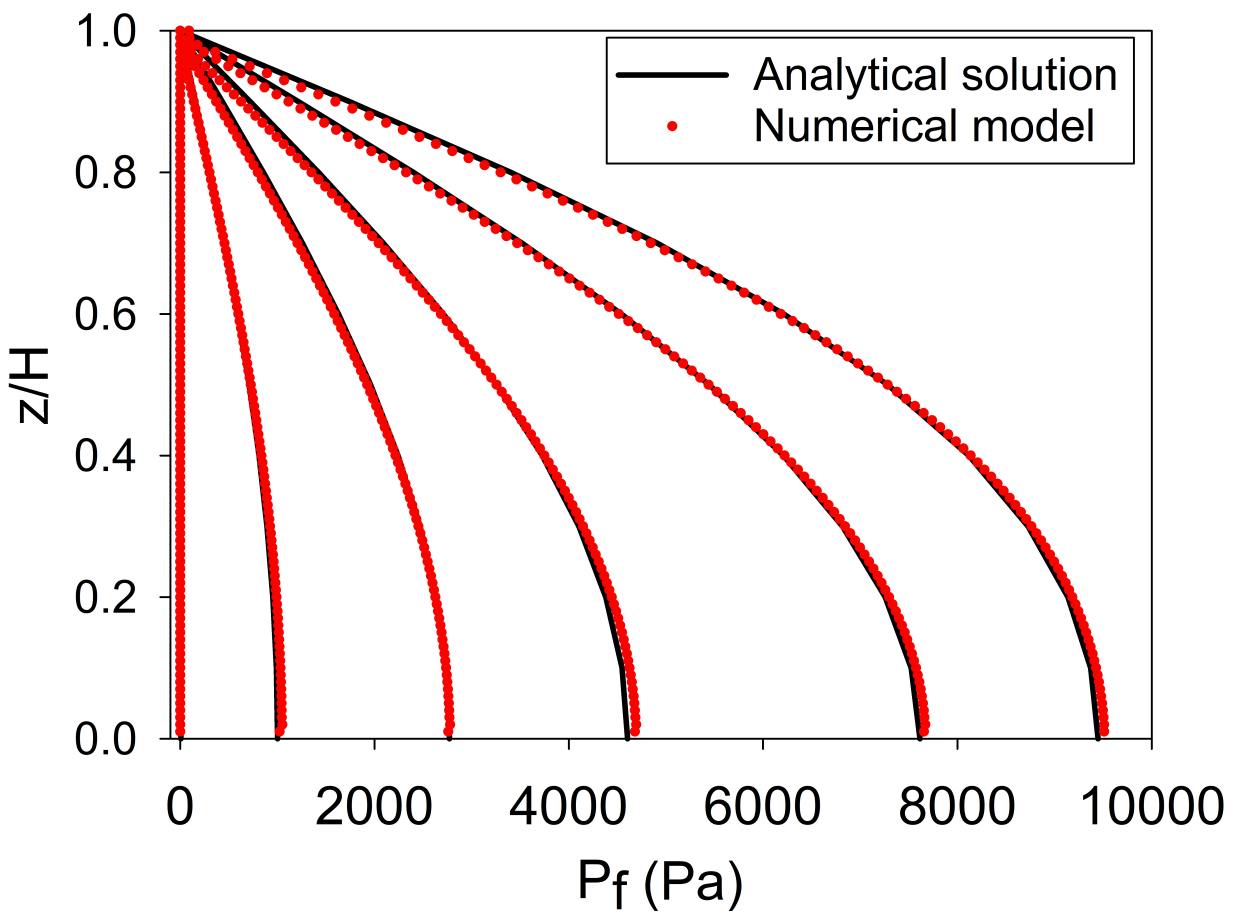

*Figure 7: Compasion between analytical solution and numerical solution*

### *Isothermal consolidation*

A common benchmark fo a fully saturated porous meida is the simulation of one-dimensional consolidation. Using the Carman-Kozeny formula, the time-dependent pressure can be caluated as:

$$p_f = \sum_{m=1}^{\infty} \frac{2F_{ext}}{M} \sin(\frac{Mz}{H}) e^{-M^2 T_V} \text{ with } M = \frac{\pi}{2}(2m+1) \tag{95}$$

where the consolidation rate $T_v = C_v t/H^2$, the consolidation coefficient $C_v = E_v n^3 d^2/(180(1-n)^2\mu)$ and the Oedometer modulus $E_v = E(1-\upsilon)/(1+\upsilon)/(1-2\upsilon)$. Our numerical model is validated by modeling the consolidation of a 1m column. This fluid has water properties (bulk modulus is 2GPa, density is 998 kg/m3 at 5 degrees Celsius and 101325 Pa (1atm) pressure, dynamic viscosity $\mu$ is 1mPa s). The porous media is modeled by elastic material with Young's modulus is 10 MPa, Poisson's ratio is 0.3, and density is 2650 kg/m3. The volume fraction of porous media $\phi_s$ is 0.7 which is equivalent to the porosity of 0.3 and the average grain diameter $d$ is 1mm. The model is discretized in 100 finite element with 1 material point per element. The external pressure applies to the top of the column is 10 kPa. Figure 7 shows a good agreement of fluid flow prediction between the theory and the model.

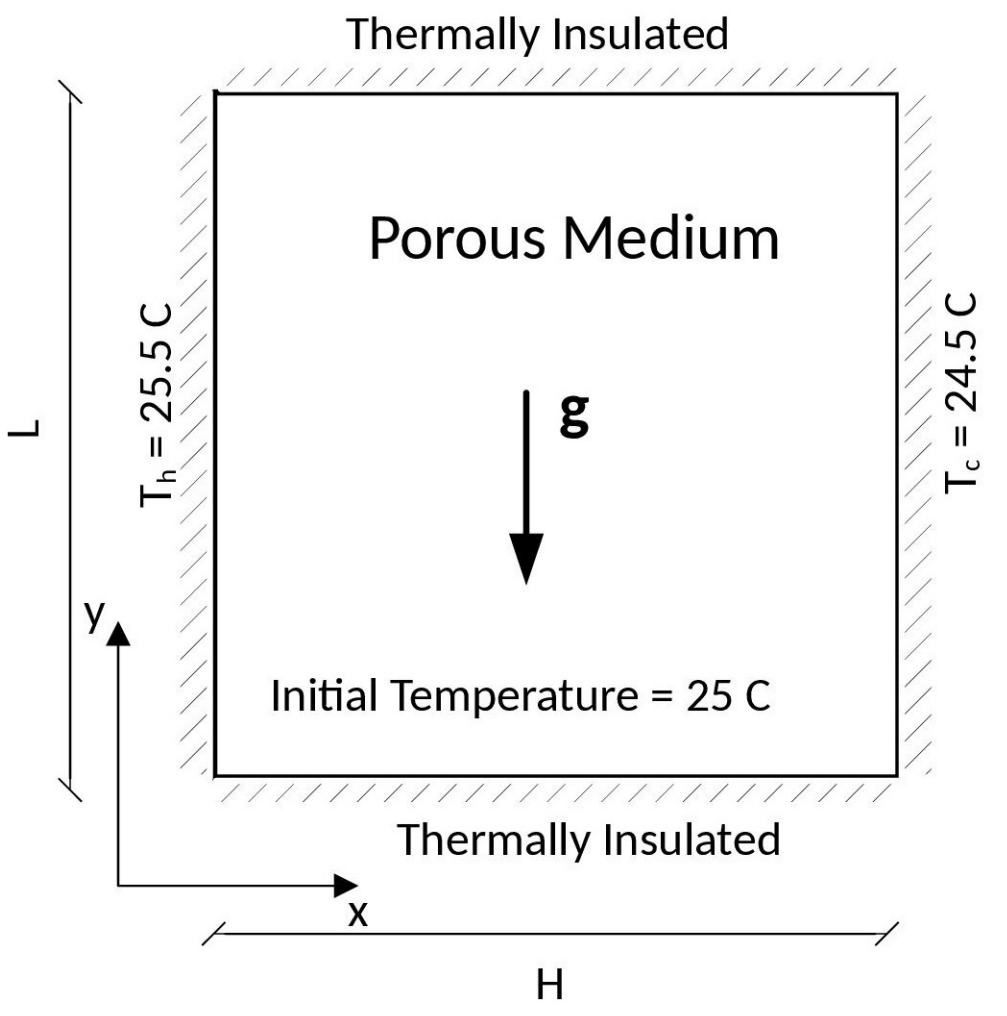

*Figure 8: Model schematic [23]*

*Thermal induced cavity flow*

Another benchkmark is the thermal induced cavity flow in porous media. Temperature and velocity distributions are calculated for a square non-deformable saturated porous media The top and bottom walls are insulated, and the left and right walls are at fixed temperatures differing by 1 degree. The fluid motion at stead state are cavity flow due to the temperature induced density variation.

The numerical is validated by comparing with the numerical solution of the finite element method. The fluid has water properties (bulk modulus is 2GPa, density is 998 kg/m3 at 5 degrees Celsius and 10325 Pa (1atm) pressure, dynamic viscosity $\mu$ is 1 mPa s). The porous media is modeled by non deformable material, and density is 2500 kg/m3. The specific heat capacity of the water and porous skeleton are 4181 J/kg.K and 835 J/kg.K respectively. The thermal conductivity of the water and porous skeleton are 0.598 W/m.K and 0.4 W/m.K. The volume fraction of porous media $\phi_s$ is 0.6 which is equivalent to the porosity of 0.4 and the average grain diameter $d$ is 1mm. The model is discretized in 20 x 20 finite element with 4 material point per element. Figure 9 shows a good agreement of numerical results of the model compared with the numerical solution of the finite element method.

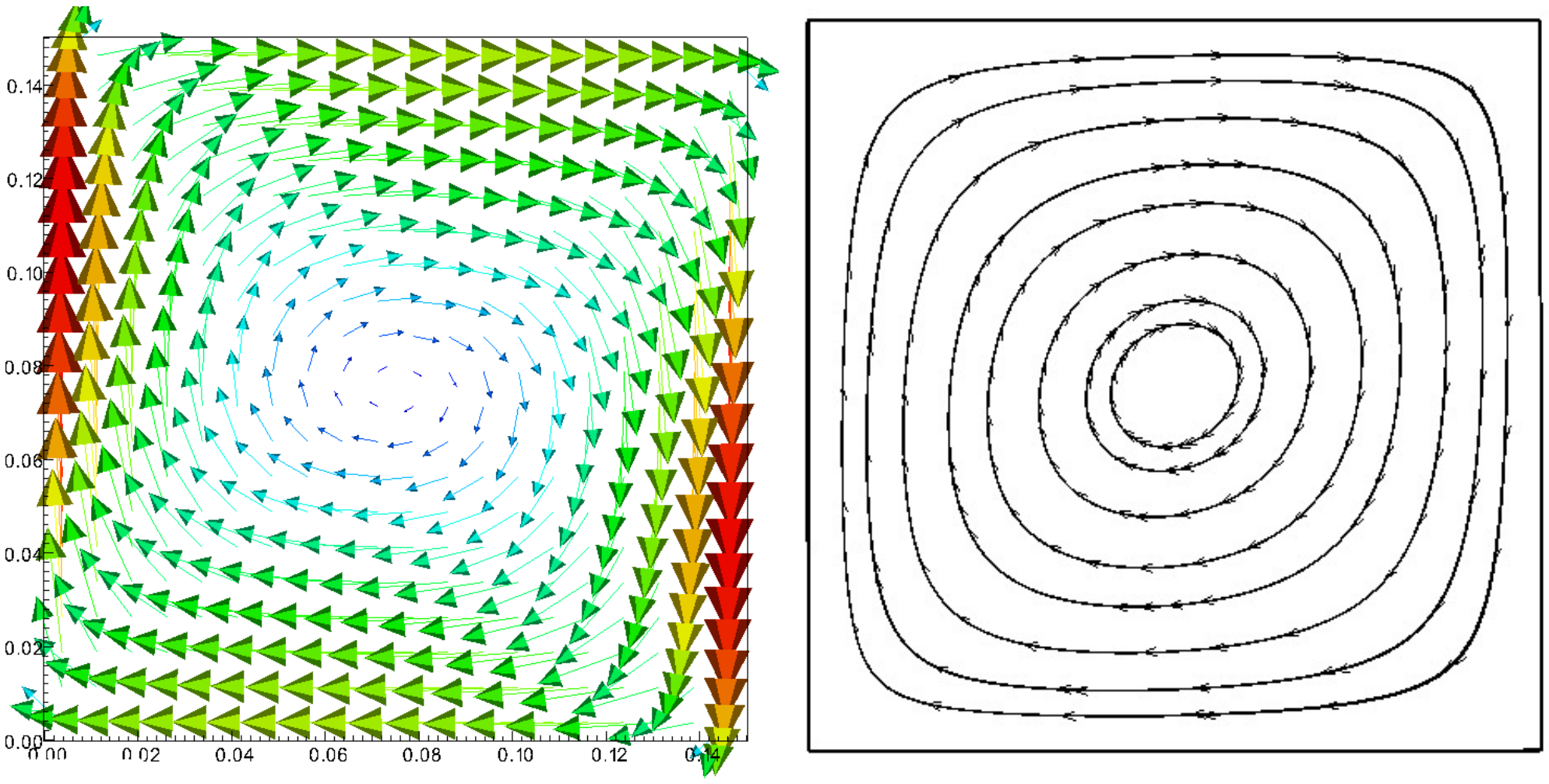

*(a) MPMICE model* *(b) FEM model [23]*

*Figure 9: Comparision between MPMICE model and FEM model*

*Underwater debris flow*

The numerical example is validated by Rzadkiewicz et al.'s experiment on submarine debris flow [24]. During the experiment, sand in a triangular box is released and then slides along a rigid bed inclined 45 degrees under water, see Figure 10.

In the numerical model, the material properties are selected based on the experiment by Rzadkiewicz et al [24]. Sand has a saturated density of 1985 $kg/m^3$ and yield stress of 200 Pa. Young's modulus has little effect on debris flow run-out because of the extreme large deformation of the debris. Therefore, we select 50 MPa Young's modulus with 0.25 Poisson's ratio. The rigid bed is much stiffer with bulk modulus and shear modulus of $117E^7$ Pa and $43.8E^7$ Pa. Under gravity, the density of the water at the surface is 999.8 $kg/m^3$ at the pressure of 1 atm. At the top boundary, the air has a density of 1.17 $kg/m^3$ at the atmospheric pressure of 1 atm. At 5 Celcius degrees, air and water have viscosity of $18.45E^{-3}$ mPa s and 1 mPa s respectively. The numerical parameters used in this example are presented in Table 1.

On all boundary faces, the Dirichlet boundary condition is imposed for velocity (u = 0 m/s) and temperature (T = 5 Celcius degrees), while the Neuman boundary condition is imposed at the top boundary for pressure (dp/dx = 0 kPa) and density (d$\rho$/dx = 0 $kg/m^3$). For the background mesh, there are

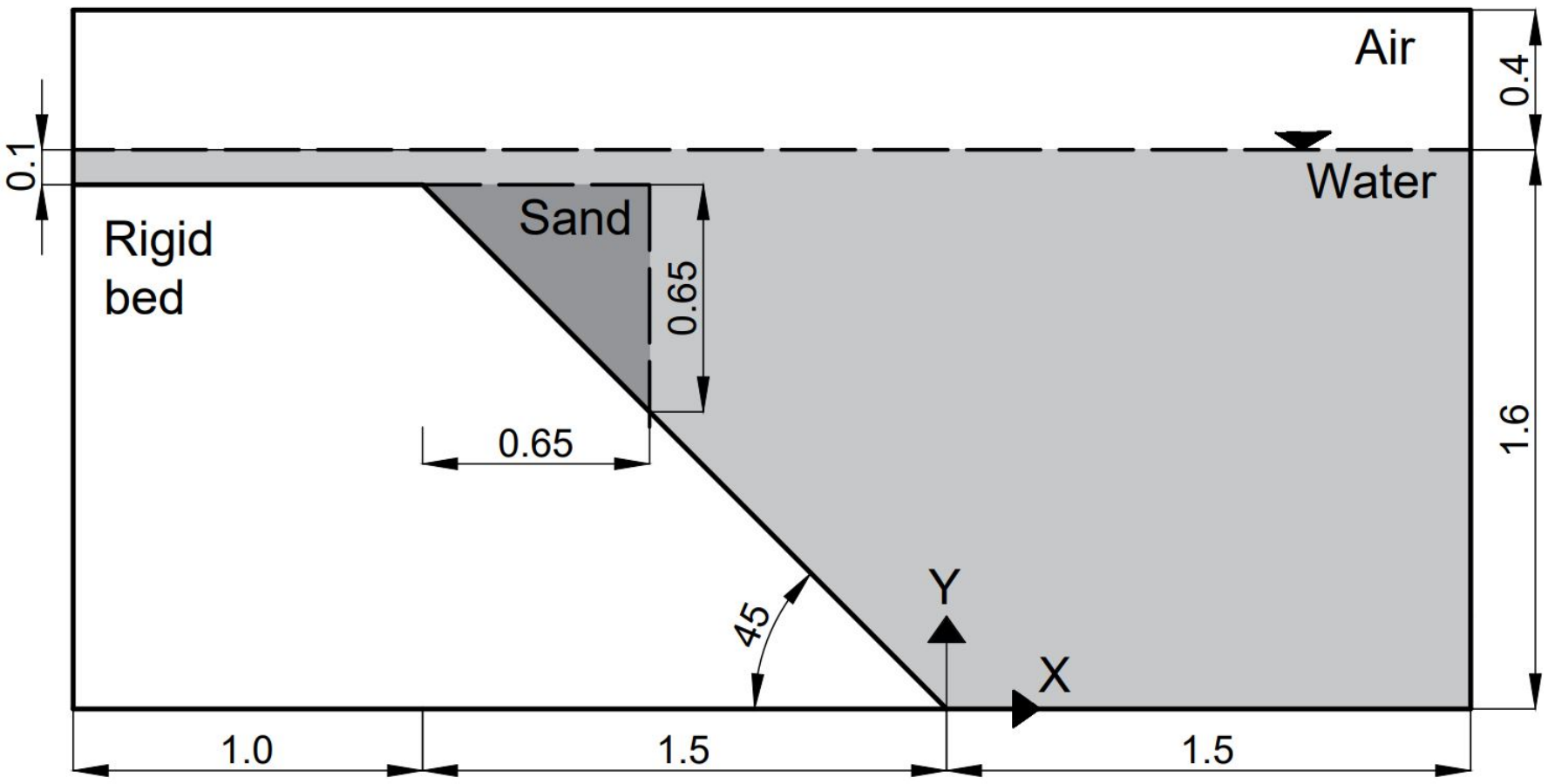

*Figure 10: Model schematic*

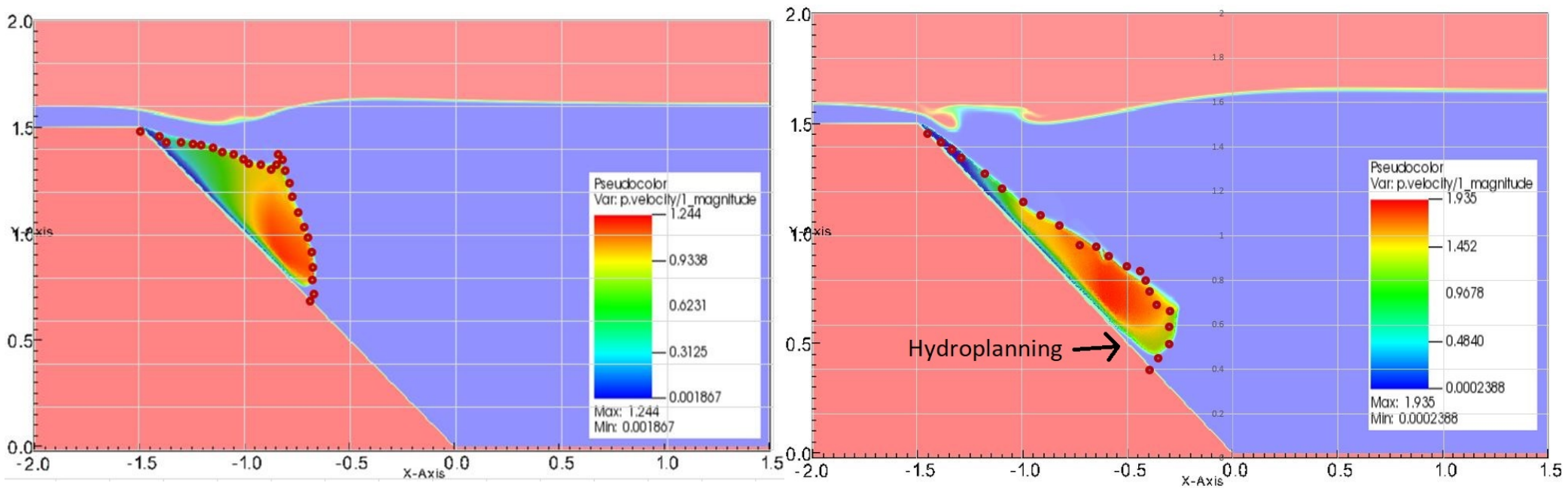

*(a) 0.4 seconds*

*(b) 0.8 seconds*

*Figure 11: Simulation of underwater debris flow*

700 x 400 = 280.000 cells. In each cell of the debris flow and rigid bed, there are 2 x 2 material points.

Figure 11a and 11b show snapshots of the debris flow sliding in the plane at 0.4 s and 0.8 s. Our simulations match the computed results from Rzadkiewicz et al. [24]. The model also captures typical hydroplaning mechanism of the underwater debris flow (hydroplaning means the debris flow is lifted up and no longer in contact with the bottom layer). The elevation of the free surface at 0.4s and 0.8s is compared between our proposed method and other methods in Figure 12. Once again, our computed results were consistent with both the experiment and others computational results [7]. Unlike other computational models based on total stress analysis, the proposed model based on the effective stress analysis which allows to analyze the water pressure and temperature in the debris flow.

| **Materials** | **Bulk modul (Pa)** | **Shear modul (Pa)** | **Density (kg/m3)** | **Temp (C)** | **Dynamic vicosity (Pa s)** | **Yield stress (Pa)** |
|---|---|---|---|---|---|---|
| Water (at surface) | 2.15e9 | - | 999.8 | 5 | 855e-6 | - |
| Air (at top boundary) | - | - | 1.177 | 5 | 18.45e-6 | - |
| Sand (porous media) | 8.33e6 | 20e6 | 1985 | 5 | - | 200 |
| Rigid bed (solid) | 117e7 | 43.8e7 | 8900 | 5 | - | - |

*Table 1: Numerical parameters for the underwater submarine debris*

We also explore the difference between underwater debris flow and saturated debris flow in terms of interacting with obstacle. Figure 13 shows the snapshot of the simulations of underwater and saturated debris flow. The saturated debris flow (see Figure 13a) behaves like frictional flow as grain have contact forces with each other. On the other hand, the underwater debris flow (see Figure 13b) behaves like tubulent flow as grains are separated from each other and exhibit no contact forces between grains.

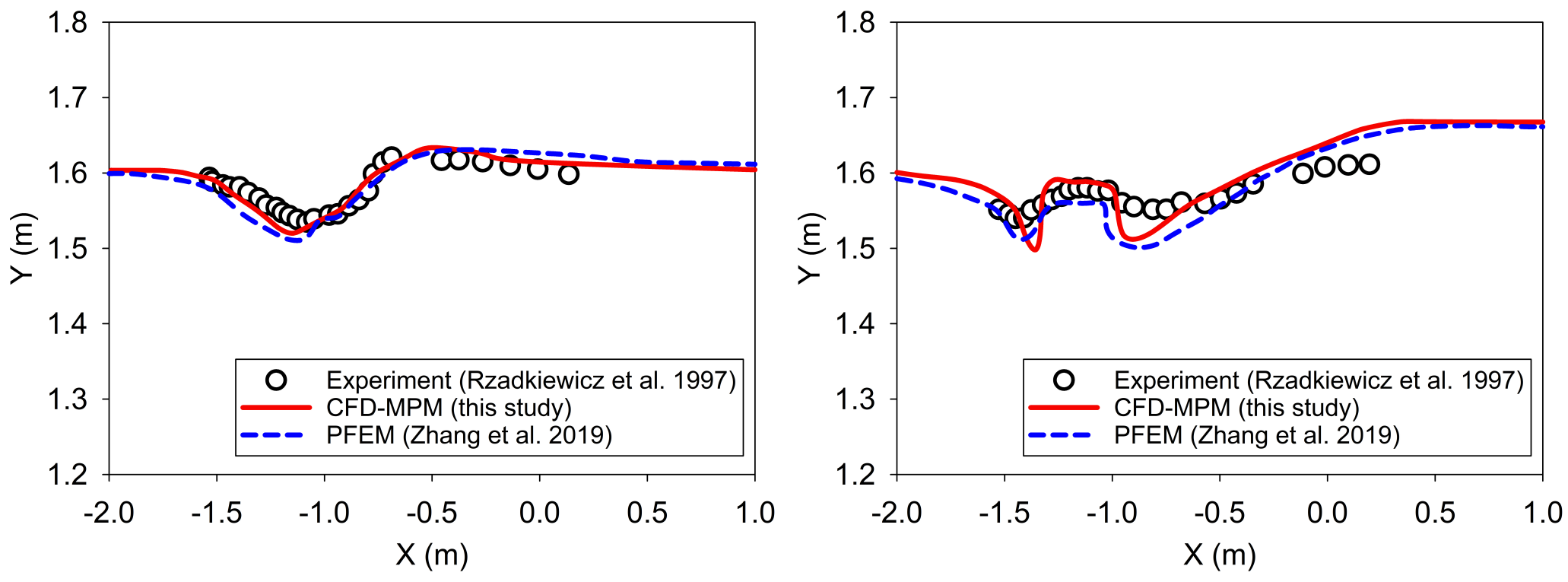

*(a) 0.4 seconds*

*(b) 0.8 seconds*

*Figure 12: Simulation of underwater debris flow*

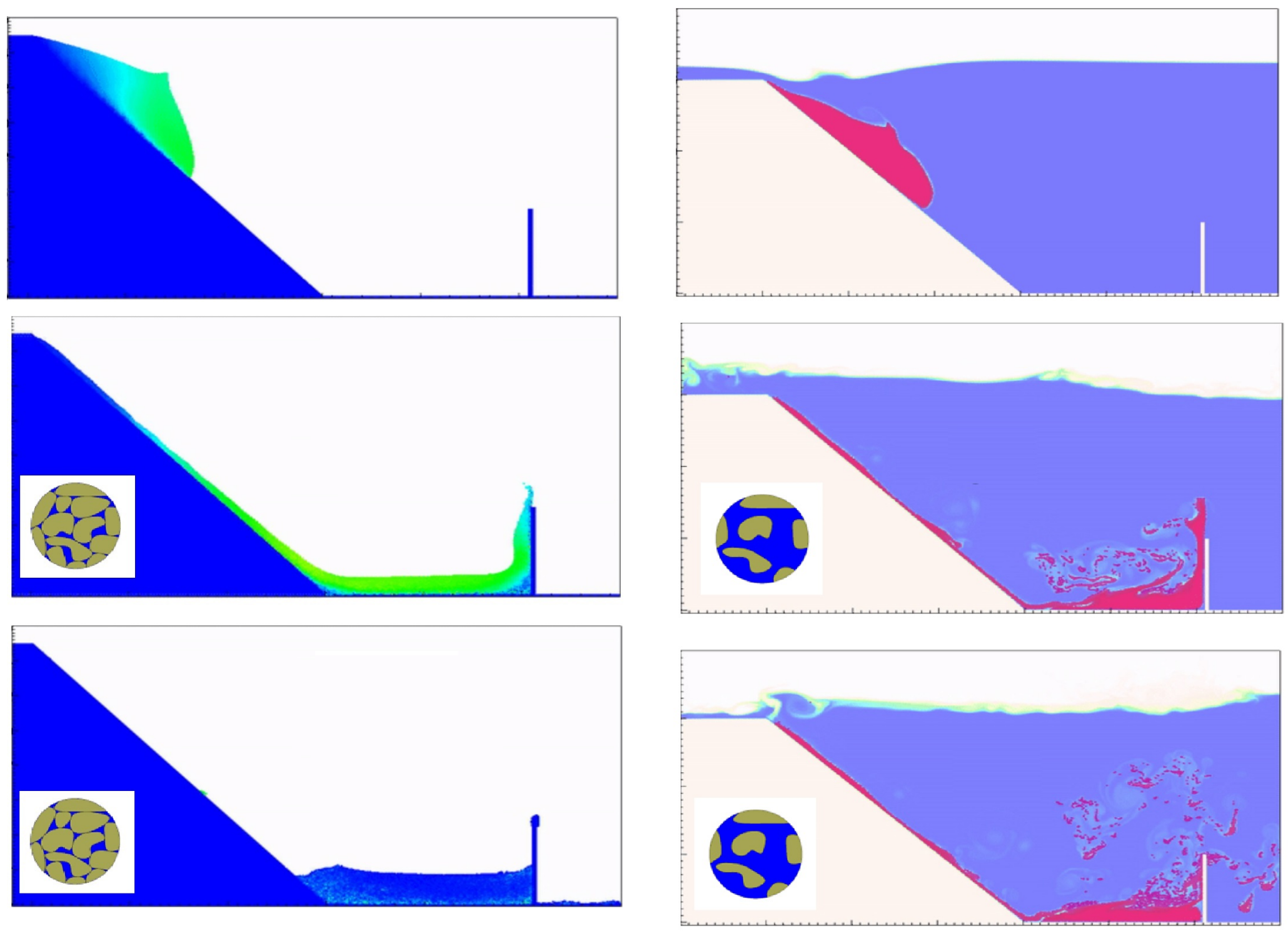

*(a) saturated debris flow using MPM*

*(b) underwater debris flow using MPMICE*

*Figure 13: Simulation of underwater debris flow*

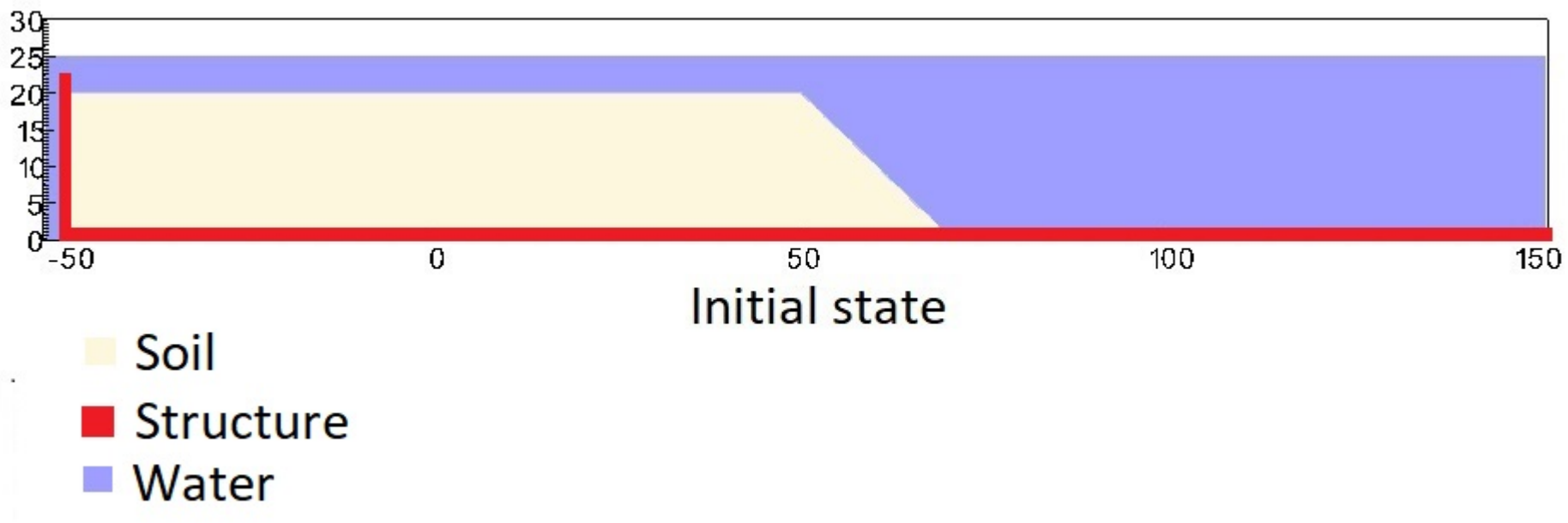

*Figure 14: Numerical model of the earthquake-induced submarine landslide*

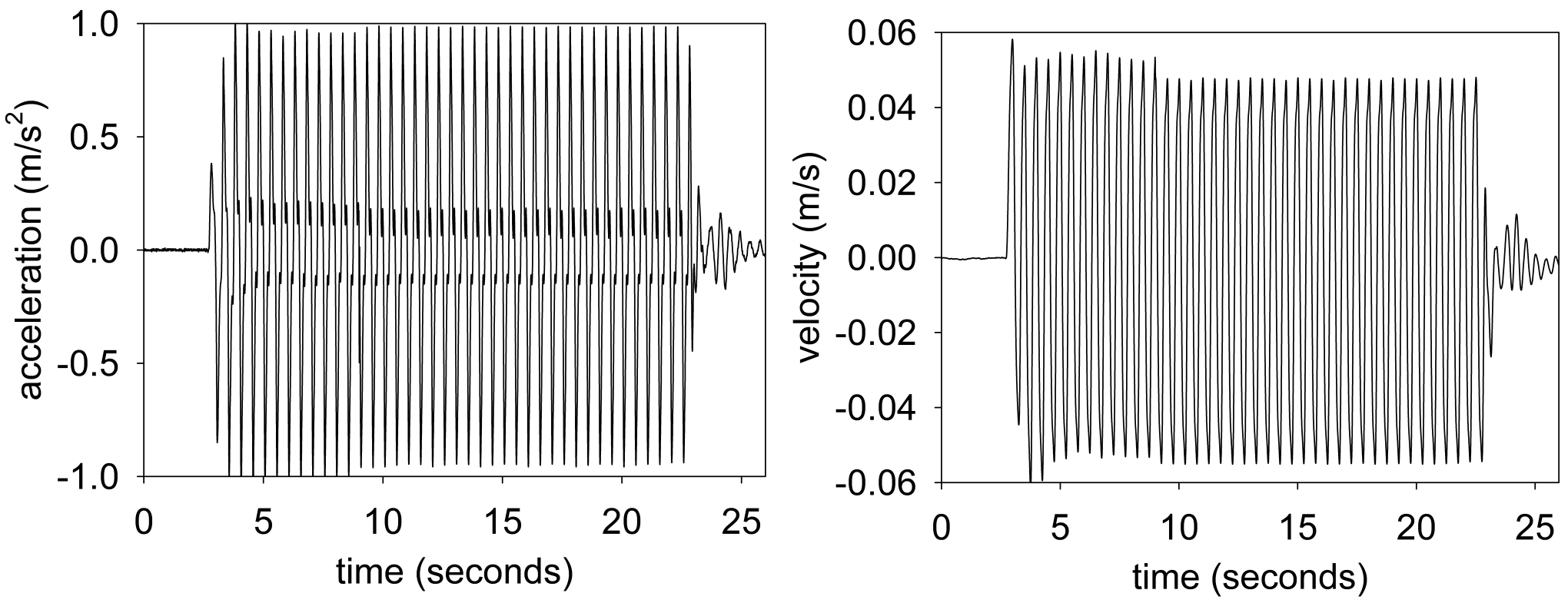

*(a) accerleration* *(b) velocity*

*Figure 15: Ground acceleration profile, frequency of 2Hz and magnitude of 1g*

In the final example, we perform numerical analysis of the earthquake induced submarine landslides. A plane strain model with the slope under water is shown in Figure 14. A 20m high slope with slope gradient of 45 degrees is placed in a horizontal and vertical structure which was used to be a skaing table to apply earthquake loading. We simplify the earthquake loading by simulating the ground shaking for 20 seconds with the peak ground acceleration of 1g and the frequency of 2Hz (Figure 15a). The ground motion is applied in terms of velocity (Figure 15b). The earthquake of this magnitude can occured typically for the earthquake of magnitude of more than 6.

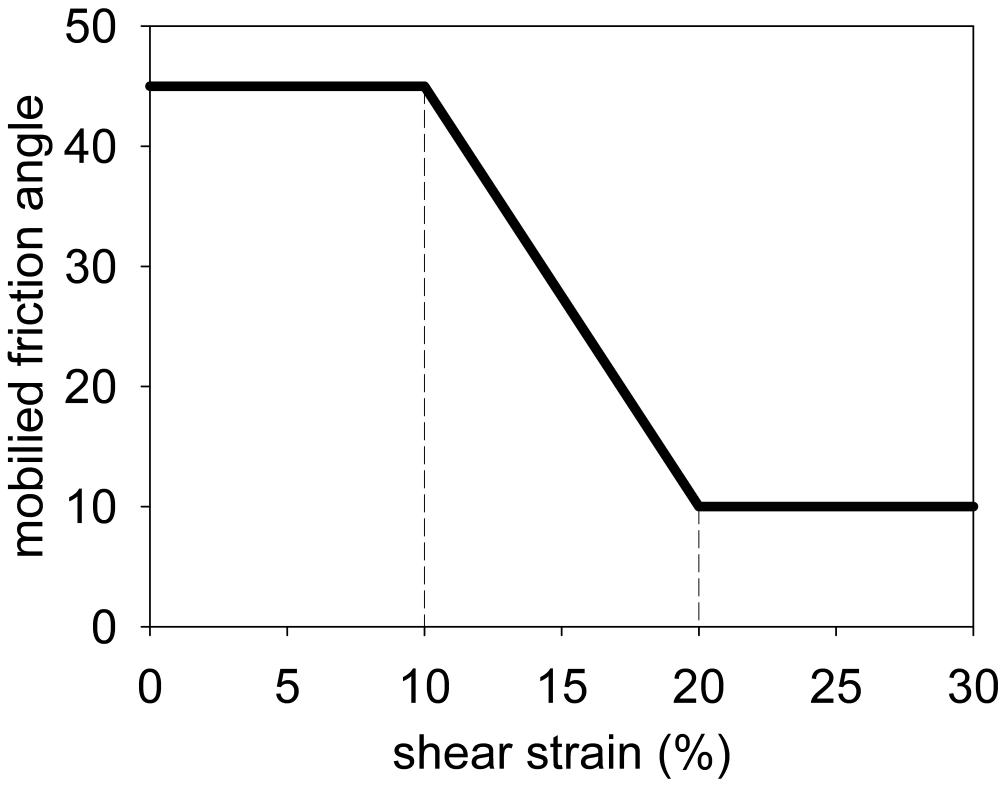

*Figure 16: Mobilized friction angle in Mohr Coulomb model*

A non-associated Mohr-Coulomb model is used for the soil. The soil grain has the density of 2650 kg/$m^3$, Young's modulus of 10 kPa and Poisson's ratio of 0.3 and zero cohesion. The mobilized friction angle $\phi'_m$ is governed following the softening curve (see Figure 16) with the peak friction angle $\phi'_p$ of 45 degrees and the residual friction angle $\phi'_r$ of 10 degrees. The porosity is 0.3 and the average grain size of the soil is around 0.1 $\mu$m to mimic the undrained behavior. The mobilized dilatancy angle is calculated from the Row-stress dilatancy as follow:

$$\sin \psi'_m = \frac{\sin \phi'_m - \sin \phi'_r}{1 - (\sin \phi'_r \sin \phi'_m)} \tag{96}$$

The solid plane is modeled as a rigid body acted as a shaking table. The contact between horizontal plane and the sand is the frictional contact with the friction coefficient of 0.1. No artificial damping is applied in the simulation. The contact between vertical plane and the sand is consdered to be smooth with zero friction coefficient. Under gravity, the density of the water at the surface is 999.8 $kg/m^3$ at the pressure of 1 atm. At the top boundary, the air has a density of 1.17 $kg/m^3$ at the atmospheric pressure of 1 atm. At 5 Celcius degrees, air and water have viscosity of $18.45e^{-3}$ mPa s and 1 mPa s respectively. On all boundary faces, the symmetric boundary condition is imposed, while the Neuman boundary condition is imposed at the top boundaryfor pressure (dp/dx = 0 kPa) and density (d$\rho$/dx = 0 $kg/m^3$). The mesh size is 0.25 x 025m with 300852 element cells and 142316 material points. The simulation takes a couple of hours to perform 60 seconds of the simulation using 4096 CPUs.

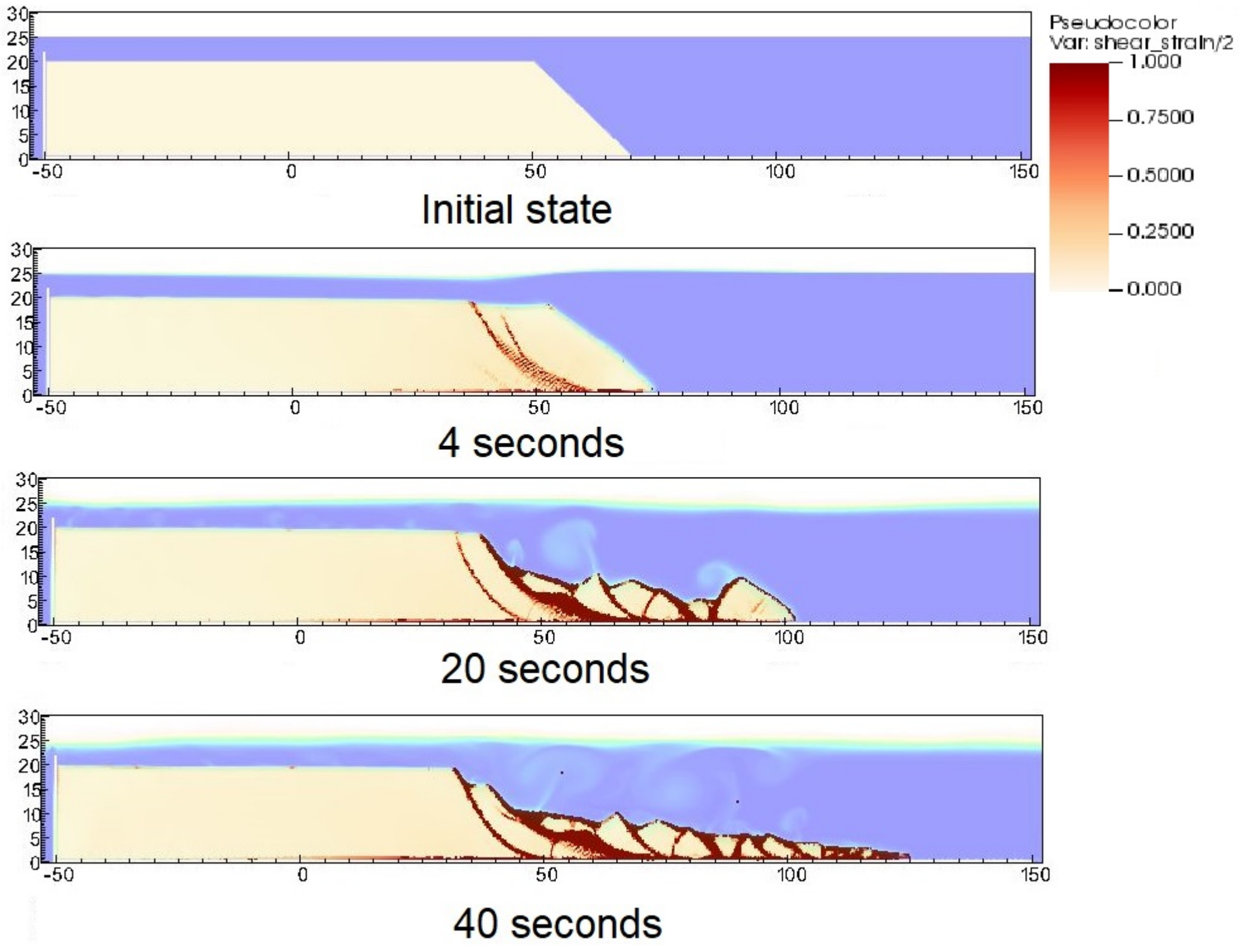

*Figure 17: Shear strain during the earthquake-induced submarine landslides*

We demonstrate the entire process and the mechanism of the earthquake-induced submarine landslides by showing the shear strain (Figure 17), the pore water ressure in atm (Figure 18) and the velocity (Figure 19). The failure mechanism can be charaterized as the progressive failure mechanism. Here are some numerical observation:

1. At the initial of the seismic event, the seismic loading triggers the first slide at 3 seconds. At 4 seconds, the debris start to move with the maximum speed of around 2-3 m/s with multiple shear band developed in the slope. The wave generated from the submarine slide is around 2-3m towards the slide direction.

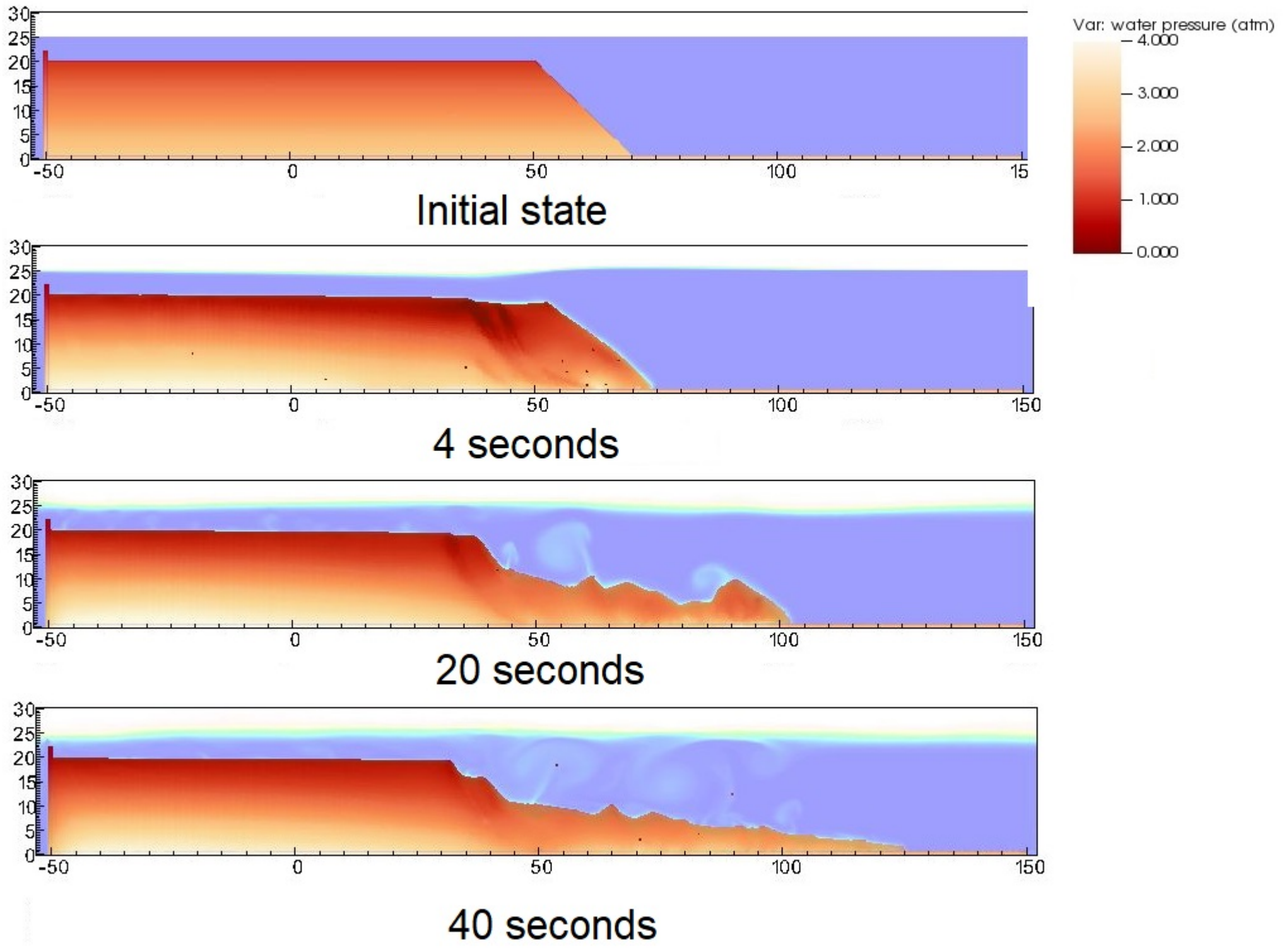

*Figure 18: pore water pressure during the earthquake-induced submarine landslides*

2. When the onset of the shear band occurs in the slope (for example at 4 seconds and 20 seconds), the negative excess pore water pressure is developed along this shear band with pore water pressure is under 1atm. This is a typical dilatancy behavior when the soil is sheared rapidly in the undrained behavior.
3. When the seismic loading ends at 23 seconds, the last shear band is mobilized and the slope soon reaches to the final deposition. No more progressive failure developed in the slope. The turbulent flow developed as the interaction between debris flow and seawater.

Overall, we show the completed process of the earthquake-induced submarine landslides involving (1) earthquake triggering mechanism, (2) the onset of the shear band with the delvelopment of negative excess pore water pressure, (3) progressive failure mechanism, (4) submarine landslide induced wave to final deposition.

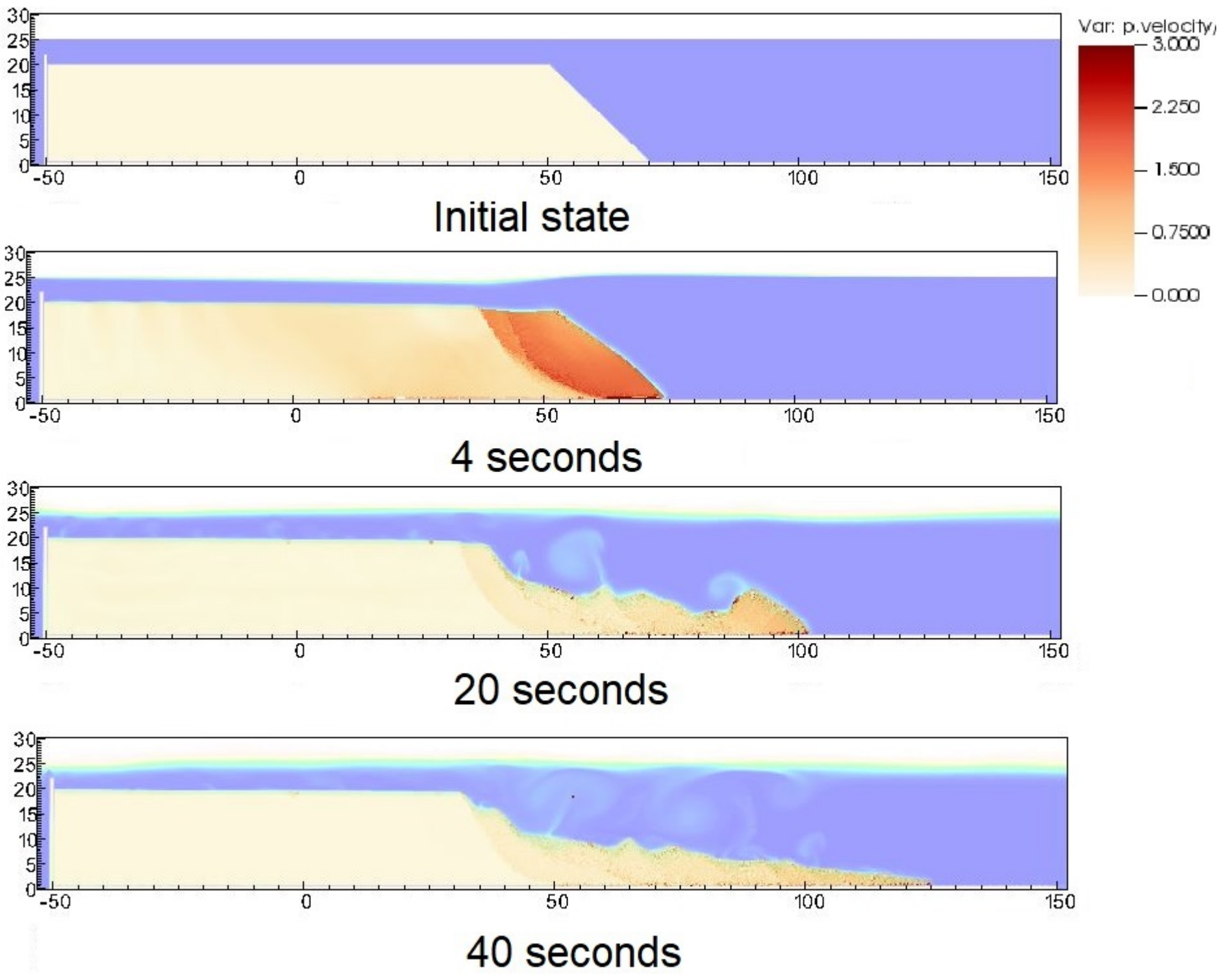

*Figure 19: Velocity during the earthquake-induced submarine landslides*

## Conclusions

We have presented a numerical approach MPMICE for the simulation of large deformation soil-fluid-structure interaction, emphasizing the simulation of the earthquake-induced submarine landslides. The model uses (1) the Material Point Method for capturing the large deformation of iso-thermal porous media and solid structures and (2) Implicit Continuous Eulerian (compressible, conservative multi-material CFD formulation) for modeling the complex fluid flow including turbulence. This model is implemented in the high-performance Uintah computational framework and validated against analytical solution and experiment. We then demonstrate the capability of the model to simulate the entire process of the earthquake induced submarine landslides.

## Acknowledgements

The authors gratefully acknowledge Dr. James Guilkey and Dr. Todd Harman from the University of Utah for sharing the insight on the theoretical framework of MPMICE which this work is based. This project has received funding from the European Union's Horizon 2020 research and innovation program under the Marie Skłodowska-Curie Actions (MSCA) Individual Fellowship (Project SUBSLIDE "Submarine landslides and their impacts on offshore infrastructures") grant agreement 101022007. The authors would also like to acknowledge the support from the Research Council of Norway through its Centers of Excellence funding scheme, project number 262644 Porelab. The computations were performed on High Performance Computing resources provided by UNINETT Sigma2 - the National Infrastructure for High Performance Computing and Data Storage in Norway.

## Appendix: Equation derivation

Before deriving the governing equation, we define the Lagrangian derivative for a state variable $f$ as:

$$\frac{D_f f}{Dt} = \frac{\partial f}{\partial t} + \boldsymbol{U_f} \cdot \nabla f \qquad \frac{D_s f}{Dt} = \frac{\partial f}{\partial t} + \boldsymbol{U_s} \cdot \nabla f \tag{97}$$

we use some definition following [16] as below:

$$-\frac{1}{V}\left[\frac{\partial V_f}{\partial p}\right] \equiv \kappa_f \qquad \text{isothermal compressibility of fluid} \tag{98}$$

$$\frac{1}{V}\left[\frac{\partial V_f}{\partial T}\right] \equiv \alpha_f \qquad \text{contant pressure thermal expansivity of fluid} \tag{99}$$

Then, the rate of volume with incompressible solid grains are calculated as below:

$$\frac{1}{V}\frac{D_f V_f}{Dt} = \frac{1}{V}\left(\left[\frac{\partial V_f}{\partial p}\right]\frac{D_f P_{eq}}{Dt} + \left[\frac{\partial V_f}{\partial T_f}\right]\frac{D_f T}{Dt}\right) = \frac{1}{V}\left(-\kappa_f \frac{D_f P_{eq}}{Dt} + \alpha_f \frac{D_f T_f}{Dt}\right) \tag{100}$$

### *Evolution of porosity*

Solving the solid mass balance equation (4) with the definition of solid mass in equation (2), it leads to the rate of porosity as below:

$$\frac{D_s m_s}{Dt} = \frac{D_s(\phi_s \rho_s V)}{Dt} = \rho_s V \frac{D_s \phi_s}{Dt} + \phi_s V \frac{D_s \rho_s}{Dt} + \phi_s \rho_s \frac{D_s V}{Dt} = 0 \tag{101}$$

The soil grains are assumed to be incompressible, therefore, term 2 in the right hand side is zero.

$$V\frac{D_s\phi_s}{Dt} + \phi_s\frac{D_sV}{Dt} = 0 \tag{102}$$

Dividing all terms with $V$ with the equation $\frac{1}{V}\frac{D_sV}{Dt} = \nabla\cdot\boldsymbol{U}_s$, it leads to:

$$\frac{D_sn}{Dt} = \frac{\partial n}{\partial t} + \boldsymbol{U}_s\cdot\nabla n = \phi_s\nabla\cdot\boldsymbol{U}_s \tag{103}$$

*Momentum conservation*

The linear momentum balance equation for the fluid phases based on mixture theory is:

$$\frac{1}{V}\frac{D_f(m_f\boldsymbol{U}_f)}{Dt} = \nabla\cdot(-\phi_f p_f\boldsymbol{I}) + \nabla\cdot\boldsymbol{\tau}_f + \overline{\rho}_f\boldsymbol{b} + \sum\boldsymbol{f}_d + \boldsymbol{f}_b \tag{104}$$

On the right hand sand, the first term is the divergence of partial fluid phase stress, the third term is the body force, the fourth term is the drag force (momentum exchange) and the fifth term is the buoyant force described in [25] for the immiscible mixtures. The buoyant force is in the form:

$$\boldsymbol{f}_b = \boldsymbol{\sigma}_f\nabla(n) \tag{105}$$

As a result, the linear momentum balance equation for the fluid phases becomes:

$$\frac{1}{V}\frac{D_f(m_f\boldsymbol{U}_f)}{Dt} = \frac{1}{V}\left[\frac{\partial(m_f\boldsymbol{U}_f)}{\partial t} + \nabla\cdot(m_f\boldsymbol{U}_f\boldsymbol{U}_f)\right] = -\phi_f\nabla p_f + \nabla\cdot\boldsymbol{\tau}_f + \overline{\rho}_f\boldsymbol{b} + \sum\boldsymbol{f}_d \tag{106}$$

The Reynolds stress component can be included in the term $\boldsymbol{\tau}_f$ to consider the turbulent effects if needed. To derive the linear momentum balance equation for the solid phase, we begin with the linear momentum balance equation for the mixture as:

$$\frac{1}{V}\frac{D_f(m_f\boldsymbol{U}_f)}{Dt} + \frac{1}{V}\frac{D_s(m_s\boldsymbol{U}_s)}{Dt} = \nabla\cdot(\boldsymbol{\sigma}) + \overline{\rho}_f\boldsymbol{b} + \overline{\rho}_s\boldsymbol{b} \tag{107}$$

Combining Terzaghi's equation (3) and subtracting both sides with equation (106), we obtain the linear momentum balance equation for the solid phase as:

$$\frac{1}{V}\frac{D_s(m_s\boldsymbol{U}_s)}{Dt} = \nabla\cdot(\boldsymbol{\sigma}') - \phi_s\nabla p_f + \overline{\rho}_s\boldsymbol{b} - \sum\boldsymbol{f}_d + \sum\boldsymbol{f}_{fric} \tag{108}$$

Here the $\boldsymbol{f}_{fric}$ stems from the soil-structure interaction following the contact law between the soil/structure interaces.

*Energy conservation*

We adopt the general form of the total energy balance equation for the porous media from [26], the total energy balance equations for the fluid phases are:

$$\frac{1}{V}\frac{D_f(m_f(e_f+0.5\boldsymbol{U}_f^2))}{Dt} = \nabla\cdot(-np_f\boldsymbol{I})\cdot\boldsymbol{U}_f+\nabla\cdot\boldsymbol{q}_f+(\overline{\rho}_f\boldsymbol{b})\cdot\boldsymbol{U}_f+\sum\boldsymbol{f}_d\cdot\boldsymbol{U}_f+\sum q_{sf} \tag{109}$$

Applying the product rule $D(m\boldsymbol{U}^2) = D(m\boldsymbol{U}\cdot\boldsymbol{U}) = 2\boldsymbol{U}\cdot D(m\boldsymbol{U})$, the left hand side of equation (109) becomes:

$$\frac{1}{V}\frac{D_f(m_f(e_f+0.5\boldsymbol{U}_f^2))}{Dt} = \frac{1}{V}\frac{D_f(m_fe_f)}{Dt}+\frac{1}{V}\frac{D_f(m_f\boldsymbol{U}_f)}{Dt}\cdot\boldsymbol{U}_f \tag{110}$$

Combining equations (106), (109), (110), we obtain the final form of the internal energy balance equation for the fluid phases as:

$$\frac{1}{V}\frac{D_f(m_fe_f)}{Dt} = \frac{1}{V}\left[\frac{\partial(m_fe_f)}{\partial t}+\nabla\cdot(m_fe_f\boldsymbol{U}_f)\right] = \overline{\rho}_fp_f\frac{D_f\upsilon_f}{Dt}+\nabla\cdot\boldsymbol{q}_f+\sum q_{sf} \tag{111}$$

On the right hand side, the terms include the average pressure-volume work, the average viscous dissipation, the thermal transport and the energy exchange between solid and fluid respectively. The heat flux is $\boldsymbol{q}_f = \overline{\rho}_f\beta_f\nabla T_f$ with $\beta_f$ being the thermal conductivity coefficient. To derive the internal energy balance equation for the solid phase, we introduce the rate of the internal energy for the thermoelastic materials as a function of elastic strain tensor $\boldsymbol{\epsilon}_s^e$ and temperature $T_s$:

$$\frac{m_s}{V}\frac{D_s(e_s)}{Dt} = \boldsymbol{\sigma}':\frac{D_s(\boldsymbol{\epsilon}_s^e)}{Dt}+\frac{D_s(e_s)}{D_s(T_s)}\frac{D_s(T_s)}{Dt} = \boldsymbol{\sigma}':\frac{D_s(\boldsymbol{\epsilon}_s^e)}{Dt}+c_v\frac{D_s(T_s)}{Dt} \tag{112}$$

$c_v$ is the specific heat at the constant volume of the solid materials. The total energy balance equation for the mixture based on [26] can be written as:

$$\begin{aligned}\frac{1}{V}\frac{D_f(m_f(e_f+0.5\boldsymbol{U}_f^2))}{Dt}+\frac{1}{V}\frac{D_s(m_s(e_s+0.5\boldsymbol{U}_s^2))}{Dt} &= \nabla\cdot(-\phi_fp_f\boldsymbol{I})\cdot\boldsymbol{U}_f \\ +\nabla\cdot(\boldsymbol{\sigma}'-\phi_sp_f\boldsymbol{I})\cdot\boldsymbol{U}_s+(-\phi_fp_f\boldsymbol{I}):\nabla\boldsymbol{U}_f+\boldsymbol{\sigma}':\nabla\boldsymbol{U}_s \\ +(\overline{\rho}_f\boldsymbol{b})\cdot\boldsymbol{U}_f+(\overline{\rho}_s\boldsymbol{b})\cdot\boldsymbol{U}_s+\nabla\cdot\boldsymbol{q}_f+\nabla\cdot\boldsymbol{q}_s+\sum\boldsymbol{f}_d\cdot(\boldsymbol{U}_f-\boldsymbol{U}_s)\end{aligned} \tag{113}$$

Subtracting equation (113), (112) to equations (109) and (108), we obtained the internal energy balance equation for solid phase as:

$$\boldsymbol{\sigma'} : \frac{D_s(\boldsymbol{\epsilon}_s^e)}{Dt} + \frac{m_s}{V} c_v \frac{D_s(T_s)}{Dt} = \Delta W_s + \Delta W_{friction} + \nabla \cdot \boldsymbol{q}_s - \sum q_{sf} \tag{114}$$

On the right hand side, the terms include the work rate from frictional sliding between solid materials $\Delta W_{friction}$, thermal transport and energy exchange between solid and fluid respectively. The heat flux is $\boldsymbol{q}_s = \overline{\rho}_s \beta_s \nabla T_s$ with $\beta_s$ being the thermal conductivity of the solid materials, the mechanical work rate $\Delta W_s = \boldsymbol{\sigma'} : \frac{D_s(\boldsymbol{\epsilon}_s)}{Dt} = \boldsymbol{\sigma'} : (\frac{D_s(\boldsymbol{\epsilon}_s^e)}{Dt} + \frac{D_s(\boldsymbol{\epsilon}_s^p)}{Dt})$ computed from the constitutive model with $\boldsymbol{\epsilon}_s^p$ is the plastic strain tensor, . By subtracting the term $\boldsymbol{\sigma'} : \frac{D_s(\boldsymbol{\epsilon}_s^e)}{Dt}$, we get the final form of the energy balance equation as:

$$\frac{m_s}{V} c_v \frac{D_s(T_s)}{Dt} = \boldsymbol{\sigma'} : \frac{D_s(\boldsymbol{\epsilon}_s^p)}{Dt} + \Delta W_{friction} + \nabla \cdot \boldsymbol{q}_s - \sum q_{sf} \tag{115}$$

*Advanced Fluid Pressure*

The discretization of the pressure equation begins with the Lagrangian face-centered velocity and the equation for the pressure

$$\overline{\rho}_{f,FC} \frac{\boldsymbol{U}_{f,FC}^{n+1} - \boldsymbol{U}_{f,FC}^{n}}{dt} = n \nabla^{FC} P_{fc}^{n+1} + \overline{\rho}_{f,FC} \boldsymbol{b} \tag{116}$$

$$\kappa \frac{dP}{dt} = \nabla^c \cdot \boldsymbol{U}_{f,FC}^{n+1} \tag{117}$$

The divergence of the equation (116) with $\nabla \cdot \boldsymbol{b} = 0$ is

$$\nabla^c \cdot \boldsymbol{U}_{f,FC}^{n+1} - \nabla^c \cdot \boldsymbol{U}_{f,FC}^{n} = \nabla^c \frac{\Delta t}{\rho_{f,FC}^n} \cdot \nabla^{FC} (P_{fc}^n + \Delta P_{fc}^n) \tag{118}$$

To solve this equation, we define the face-centered intermediate velocity $\boldsymbol{U}_{f,FC}^{*}$ as:

$$\overline{\rho}_{f,FC} \frac{\boldsymbol{U}_{f,FC}^{*} - \boldsymbol{U}_{f,FC}^{n}}{\Delta t} = n \nabla^{FC} P_{fc}^{n} + \overline{\rho}_{f,FC} \boldsymbol{b} \tag{119}$$

The divergence of the equation (119) is

$$\nabla^c \cdot \boldsymbol{U}_{f,FC}^{*} - \nabla^c \cdot \boldsymbol{U}_{f,FC}^{n} = \nabla^c \frac{\Delta t}{\rho_{f,FC}^n} \cdot \nabla^{FC} P_{fc}^n \tag{120}$$

Combining equations (117, 118, 120), it leads to

$$\left(\kappa - \nabla^c \frac{\Delta t}{\rho^n_{f,FC}} \cdot \nabla^{FC}\right) \Delta P^n_{fc} = -\nabla^c \cdot \boldsymbol{U}^*_{f,FC} \tag{121}$$

When the fluid is incompressible, $\kappa$ approaches to zero and the equation (121) becomes the Poisson's equation for the incompressible fluid flow.

*Momentum and Energy exchange with an implicit solver*

Considering the fluid momentum balance equation as

$$(m\boldsymbol{U})^{n+1}_{f,FC} = (m\boldsymbol{U})^n_{f,FC} - \Delta t(Vn\nabla^{FC}P^n_{fc} + m_f\boldsymbol{b}) + VK\Delta t(\boldsymbol{U}^{n+1}_{s,FC} - \boldsymbol{U}^{n+1}_{f,FC}) \tag{122}$$

Assuming $m^{n+1}_{f,FC} = m^n_{f,FC}$ we get

$$\boldsymbol{U}^{n+1}_{f,FC} = \boldsymbol{U}^n_{f,FC} - \Delta t(\frac{\nabla^{FC}P^n_{fc}}{\rho^n_{f,FC}} + \boldsymbol{b}) + \frac{\Delta t K}{\overline{\rho}^n_{f,FC}}(\boldsymbol{U}^{n+1}_{s,FC} - \boldsymbol{U}^{n+1}_{f,FC}) \tag{123}$$

As defined in the section 'Advanced Fluid Pressure', the face-centered intermediate fluid velocity $\boldsymbol{U}^*_{f,FC} = \Delta t(\nabla^{FC}P^n_{fc}/\rho^n_{f,FC} + \boldsymbol{b})$ leading to

$$\boldsymbol{U}^{n+1}_{f,FC} = \boldsymbol{U}^*_{f,FC} + \frac{\Delta t K}{\overline{\rho}^n_{f,FC}}(\boldsymbol{U}^{n+1}_{s,FC} - \boldsymbol{U}^{n+1}_{f,FC}) \tag{124}$$

Considering the solid momentum balance equation as

$$(m\boldsymbol{U})^{n+1}_{s,FC} = (m\boldsymbol{U})^n_{s,FC} - \Delta t(V\nabla^{FC}\cdot\boldsymbol{\sigma}'^n - V(1-n)\nabla^{FC}P^n_{fc} + m_s\boldsymbol{b}) - VK\Delta t(\boldsymbol{U}^{n+1}_{s,FC} - \boldsymbol{U}^{n+1}_{f,FC}) \tag{125}$$

We define the face-centered intermediate solid velocity as $\boldsymbol{U}^*_{s,FC} = \Delta t(\nabla^{FC} \cdot \boldsymbol{\sigma}'^n_c/\overline{\rho}_{s,FC} - \nabla^{FC}P^n_{fc}/\rho_s + \boldsymbol{b})$ leading to

$$\boldsymbol{U}^{n+1}_{s,FC} = \boldsymbol{U}^*_{s,FC} - \frac{\Delta t K}{\overline{\rho}^n_{s,FC}}(\boldsymbol{U}^{n+1}_{s,FC} - \boldsymbol{U}^{n+1}_{f,FC}) \tag{126}$$

Combining equation (124) and (126) we get

$$\begin{aligned}\boldsymbol{U}^*_{f,FC} + \Delta\boldsymbol{U}_{f,FC} &= \boldsymbol{U}^*_{f,FC} + \frac{\Delta t K}{\overline{\rho}^n_{f,FC}}(\boldsymbol{U}^*_{s,FC} + \Delta\boldsymbol{U}_{s,FC} - \boldsymbol{U}^*_{f,FC} - \Delta\boldsymbol{U}_{f,FC}) \\ \boldsymbol{U}^*_{s,FC} + \Delta\boldsymbol{U}_{s,FC} &= \boldsymbol{U}^*_{s,FC} - \frac{\Delta t K}{\overline{\rho}^n_{s,FC}}(\boldsymbol{U}^*_{s,FC} + \Delta\boldsymbol{U}_{s,FC} - \boldsymbol{U}^*_{f,FC} - \Delta\boldsymbol{U}_{f,FC})\end{aligned} \tag{127}$$

Rearranging the equation (127), it leads to the linear system of equations

$$\begin{vmatrix} (1+\beta_{12,FC}) & -\beta_{12,FC} \\ -\beta_{21,FC} & (1+\beta_{21,FC}) \end{vmatrix} \begin{vmatrix} \Delta \boldsymbol{U}_{f,FC} \\ \Delta \boldsymbol{U}_{s,FC} \end{vmatrix} = \begin{vmatrix} \beta_{12,FC}(\boldsymbol{U}^*_{s,FC} - \boldsymbol{U}^*_{f,FC}) \\ \beta_{21,FC}(\boldsymbol{U}^*_{f,FC} - \boldsymbol{U}^*_{s,FC}) \end{vmatrix}$$

Solving this linear equations with $\beta_{12,FC} = (\Delta t K)/\overline{\rho}^n_{f,FC}$ and $\beta_{21,FC} = (\Delta t K)/\overline{\rho}^n_{s,FC}$ with K is the momentum exchange coefficient. Similar derivation can be performed to computed the cell-center velocity increment leading to

$$\begin{vmatrix} (1+\beta_{12c}) & -\beta_{12c} \\ -\beta_{21c} & (1+\beta_{21c}) \end{vmatrix} \begin{vmatrix} \Delta \boldsymbol{U}_{fc} \\ \Delta \boldsymbol{U}_{sc} \end{vmatrix} = \begin{vmatrix} \beta_{12c}(\boldsymbol{U}^*_{sc} - \boldsymbol{U}^*_{fc}) \\ \beta_{21c}(\boldsymbol{U}^*_{fc} - \boldsymbol{U}^*_{sc}) \end{vmatrix}$$

with $\beta_{12c} = (\Delta t K)/\overline{\rho}^n_{fc}$ and $\beta_{21c} = (\Delta t K)/\overline{\rho}^n_{sc}$ and the cell-centered intermediate velocity can be calculated by

$$\begin{aligned} \boldsymbol{U}^*_{fc} &= \boldsymbol{U}^n_{fc} + \Delta t\left(-\frac{\nabla P^{n+1}_{fc}}{\rho^n_{fc}} + \frac{\nabla \cdot \boldsymbol{\tau}^n_{fc}}{\overline{\rho}^n_{fc}} + \boldsymbol{b}\right) \\ \boldsymbol{U}^*_{sc} &= \boldsymbol{U}^n_{sc} + \Delta t\left(\frac{\nabla \cdot \boldsymbol{\sigma}'^n_c}{\overline{\rho}^n_{sc}} - \frac{\nabla P^{n+1}_{fc}}{\rho_s} + \boldsymbol{b}\right) \end{aligned} \tag{128}$$

For generalize multi materials i,j = 1:N, the linear equations is in the form:

$$\begin{vmatrix} (1+\beta_{ij}) & -\beta_{ij} \\ -\beta_{ji} & (1+\beta_{ji}) \end{vmatrix} \begin{vmatrix} \Delta \boldsymbol{U}_{i} \\ \Delta \boldsymbol{U}_{j} \end{vmatrix} = \begin{vmatrix} \beta_{ij}(\boldsymbol{U}^*_{i} - \boldsymbol{U}^*_{j}) \\ \beta_{ji}(\boldsymbol{U}^*_{j} - \boldsymbol{U}^*_{i}) \end{vmatrix}$$

Similar approach applied for the ernergy exchange term leading to:

$$\begin{vmatrix} (1+\eta_{ij}) & -\eta_{ij} \\ -\eta_{ji} & (1+\eta_{ji}) \end{vmatrix} \begin{vmatrix} \Delta T_{i} \\ \Delta T_{j} \end{vmatrix} = \begin{vmatrix} \eta_{ij}(T^n_{i} - T^n_{j}) \\ \eta_{ji}(T^n_{j} - T^n_{i}) \end{vmatrix}$$

with $\eta$ is the energy exchange coefficient.